\documentclass[10pt]{article}
\usepackage{amsthm,amsmath,amssymb,amsxtra}
\usepackage{mathrsfs}
\usepackage{float,verbatim}
\usepackage{threeparttable,booktabs}
\usepackage{graphicx}
\usepackage{epstopdf}
\usepackage{mathtools}
\usepackage{subfig}
\usepackage{algorithm,algorithmic}
\usepackage{color}
\graphicspath{{./pics/}}
\usepackage{geometry}
\usepackage{array}
\newtheorem{theorem}{Theorem}[section]
\newtheorem{lemma}{Lemma}[section]
\newtheorem{remark}{Remark}[section]
\newtheorem{ass}{Assumption}[section]

\usepackage{graphicx}
\usepackage{dsfont}
\usepackage{epstopdf}
\graphicspath{{./pics/}}
\usepackage{booktabs}
\usepackage{threeparttable}
\usepackage{verbatim}
\newcommand{\bel}{\begin{equation} \label}
\newcommand{\ee}{\end{equation}}
\def\beq{\begin{equation}}
\def\eeq{\end{equation}}
\newcommand{\bea}{\begin{eqnarray}}
\newcommand{\eea}{\end{eqnarray}}
\newcommand{\beas}{\begin{eqnarray*}}
\newcommand{\eeas}{\end{eqnarray*}}

\newcommand{\p}{\partial}

\allowdisplaybreaks
\renewcommand{\d}{{\rm d}}
\newcommand{\Om}{\Omega}

\numberwithin{equation}{section}

\renewcommand{\d}{\,\mathrm{d}}
\allowdisplaybreaks

\def\phi {\varphi}

\allowdisplaybreaks

\allowdisplaybreaks

\title{Unique and Stable Recovery of Space-Variable Order in Multidimensional Subdiffusion\thanks{The work of B. Jin is supported by Hong Kong RGC General Research Fund (Projects 14306423 and 14306824), ANR / RGC Joint Research
Scheme (A-CUHK402/24) and a start-up fund from The Chinese University of Hong Kong. The work of Y. Kian is supported by the French National Research Agency ANR and Hong Kong RGC Joint Research Scheme for the project IdiAnoDiff (grant ANR-24-CE40-7039).}}

\author{Jiho Hong\thanks{Department of Mathematics, The Chinese University of Hong Kong, Shatin, New Territories, Hong Kong SAR, P.R. China (jihohong@cuhk.edu.hk, b.jin@cuhk.edu.hk).}\and Bangti Jin\footnotemark[2]\and Yavar Kian\thanks{Univ Rouen Normandie, CNRS, Normandie Univ, LMRS UMR 6085, F-76000 Rouen, France (\texttt{yavar.kian@univ-rouen.fr})}}

\date{\today}

\begin{document}

\maketitle

\begin{abstract}
In this work we investigate the unique identifiability and stable recovery of a spatially dependent variable-order in the subdiffusion model from the boundary flux measurement. We establish several new unique identifiability results from the observation at one point on the boundary without / with the knowledge of medium properties, and a conditional Lipschitz stability estimate when the observation is available on the whole boundary. The analysis crucially employs resolvent estimates in the $L^r(\Omega)$ ($r>d$) spaces, solution representation in the Laplace domain and novel asymptotic expansions of the Laplace transform of the boundary flux at $p= 0$ and $p=1$.
\end{abstract}

\section{Introduction}
Let $\Om\subset \mathbb{R}^d$ ($d=2,3$) be an open bounded and connected domain with a $C^{1,1}$ boundary. Let
$\alpha:\Omega \to(0,1)$ be a spatially variable order function. Consider the following subdiffusion model:
\begin{equation}\label{eq:2D3D:ibvp}
\left\{
\begin{aligned}
\rho(x)\partial_t^{\alpha(x)}U(t,x) - \nabla_x\cdot( \sigma(x) \nabla_x U(t,x)) + q(x) U(t,x)&=0,\quad (t,x)\in(0,\infty)\times \Omega,\\
U(t,x)&=g(t,x),\quad (t,x)\in(0,\infty)\times\partial\Omega,\\
U(0,x)&=0,\quad x\in\Omega,
\end{aligned}\right.
\end{equation}
where $\rho\in L^\infty(\Omega)$ and $\sigma\in C^{1}(\overline{\Omega})$ have positive lower bounds, $q\in L^\infty(\Omega)$ is nonnegative, and $g$ is the boundary excitation. In the model \eqref{eq:2D3D:ibvp}, the spatially variable order Caputo fractional derivative $\partial_t^{\alpha(x)} U(t,x)$ in time $t$ is defined by (see, e.g., \cite[p. 92]{KilbasSrivastavaTrujillo:2006} or \cite[p. 41]{Jin:2021})
\begin{equation}
    \partial_t^{\alpha(x)}U(t,x) = \frac{1}{\Gamma(1-\alpha(x))}\int_0^t{(t-s)^{-\alpha(x)}}\partial_sU(s,x)\,{\rm d}s,
\end{equation}
where $\Gamma(z)$, for $\Re(z)>0$, denotes Euler's Gamma function.

The model \eqref{eq:2D3D:ibvp} describes space-dependent anomalous diffusion processes in complex media in which  heterogeneous regions exhibit spatially inhomogeneous variations. Chechkin et al \cite{ChechkinGorenflo:2005}  first studied a composite system with two regions; see also \cite{Fedotov:2012,KorabelBarkai:2010,Stickler:2011} for related studies in physics. So far such models have been successfully used in, e.g., reinforced concrete, liquid infiltration in porous
media, and nonlinear contact phenomena; see the reviews \cite{Patnaik:2020,SunChangChen:2019} for details. See also the works
\cite{Orsingher:2018,ZhangLiLuo:2013} for the derivation of variable-order models in the framework of continuous time random walk (with a space-dependent diffusion coefficient).

The concerned inverse problem is to identify the variable order $\alpha$ in the model \eqref{eq:2D3D:ibvp} from the over-posed boundary flux data. In this work, we establish several new uniqueness and stability results for the inverse problem. First, we prove that with only one point measurement and a piecewise constant variable order $\alpha(x)$, the range of $\alpha(x)$ is uniquely determined even if the medium properties (e.g., $\rho$, $\sigma$ and $q$) are unknown, and also the interfaces of discontinuity of  $\alpha(x)$ are uniquely determined when the medium properties are known, under no additional assumptions on the interfaces; see Theorem \ref{theorem:flux:lowfreq:asympt:onept} for the precise statement. Second, in a generalized setting where the variable order $\alpha$ is not necessarily piecewise constant, we prove the unique identifiability in the same sense under a monotonicity condition when the medium property is known. See Theorem \ref{theorem:uniqueness:measfuns} for the precise statement. Third, we establish a semi-continuity and conditional stability result from the flux measurement at one point on the boundary $\partial\Omega$ and  over the whole boundary $\partial\Omega$ in Theorems \ref{theorem:stability} and \ref{theorem:stability:fullbdydata}, respectively.

The analysis is based on asymptotic expansions of the Laplace transform $\p_\nu \widehat{U}(p,x)$ of the normal derivative $\partial_\nu U(t,x)$ of the solution $U(t,x)$ to problem \eqref{eq:2D3D:ibvp} as $p\to 0^+$ (for the uniqueness) and $p$ close to $1$ (for the continuity and stability) given in Theorems \ref{theorem:flux:lowfreq:asympt:general} and \ref{theorem:flux:freqisone:asympt:general}, respectively. The choice of the limiting values $p=0$ and $p=1$ is based on the key observation that for $p=0$ and $p=1$, the value of $p^{\alpha}$ is independent of $\alpha$.
The proof of the asymptotics employs suitable $L^r(\Omega)$ resolvent estimates for the following Schr\"odinger equation
\begin{equation*}
\left\{\begin{aligned}
- \nabla \cdot( \sigma \nabla\widehat{U}(p,\cdot)) + (q + \rho p^{\alpha(x)})\widehat{U}(p,\cdot)&=0,\quad\mbox{in } \Omega,\\
\widehat{U}(p,\cdot)&=\widehat{g}(p,\cdot),\quad\mbox{on }\partial\Omega,
\end{aligned}\right.
\end{equation*}
for $p>0$, which arises from the Laplace transform of problem \eqref{eq:2D3D:ibvp}, cf. \eqref{eq:afterLaptrans} below. The asymptotic expansion along with the strong maximum principle and Hopf's lemma for elliptic problems yields the uniqueness, and together with an integral representation of the boundary flux data, leads also to the continuity and stability results. Note that the potential $q(x) + \rho p^{\alpha(x)}$ in the Schr\"{o}dinger operator has infinitely many values depending on the frequency $p>0$, and thus one-point flux measurement for the model \eqref{eq:2D3D:ibvp} is more informative about the variable order $\alpha(x)$ than the one-point measurement for the standard elliptic model.

Now we situate this study in the context of inverse problems for subdiffusion, which has received a lot of attention in recent years; see the reviews \cite{JinRundell:2015,LiLiuYamamoto:2019,LiYamamoto:2019} and the monograph \cite{KaltenbacherRundell:2023} for details. However, so far the study on the variable-order subdiffusion is very limited \cite{HJK:2025:ISDVO,IkehataKian:2023,Kian:2018:TFD}, due to the unprecedented analytical challenges in the analysis of the direct problem. Kian, Soccorsi and Yamamoto \cite{Kian:2018:TFD} established the well-posedness of the model using Laplace transform, and proved the unique recovery of the variable order $\alpha(x)$ from the partial Dirichlet to Neumann map (along with two other coefficients). Ikehata and Kian \cite{IkehataKian:2023} proved the unique recovery of the geometry of the variable order from the Neumann data when the problem is equipped with a specially designed Dirichlet datum. See also the work \cite{ZhengChengWang:2019} for related results on models involving a time-dependent variable order. Very recently, in the one-dimensional case, the authors \cite{HJK:2025:ISDVO}  proved two unique identifiability results for piecewise constant variable orders from one single-point boundary flux observation using the idea of gluing. In this work, we have substantially extended the results in \cite{Kian:2018:TFD} and \cite{HJK:2025:ISDVO}: using only one point measurement, we prove novel uniqueness and moreover stability results in the multi-dimensional case. The analysis is based on resolvent estimates in $L^r(\Omega)$, $r>d$, and the asymptotics of the Laplace transform $\partial_\nu \widehat{U}(p,x)$ in time at frequency $p$ of the flux as $p\to0^+$ and $p\to 1$, respectively. The latter expansion also enables us to establish the continuity and stability results. To the best of our knowledge, there are no known continuity and stability results for the recovery of the variable order in the existing literature.

The rest of the paper is organized as follows. In Section \ref{sec:main}, we state the main results of the work, and provide further discussions on related results. In Section \ref{sec:direct}, we develop the main tools for the analysis, i.e., asymptotic expansion of the Laplace transform of the Neumann boundary data. In Section \ref{sec:inverse}, we prove the uniqueness and stability results. In Section \ref{sec:concl} we conclude the paper with further discussions. Throughout, we use the standard notation for Sobolev spaces $W^{s,r}(\Omega)$ for $s\geq0$ and $r\geq1$ \cite{Adams:2003:SS,Grisvard:1985:EPND}. We denote by $\nu$ the unit outward normal vector to the boundary $\partial\Omega$ and $\p_\nu$ the corresponding normal derivative.

\section{Main results and discussions}\label{sec:main}

In this section we describe the main results and provide further discussions on related  results in existing studies. First we state the assumptions on the problem data.

\begin{ass}\label{ass:coef}The boundary $\partial\Omega$ is of class $C^{1,1}$,
$\rho\in L^\infty(\Omega)$ and $\sigma\in C^{1}(\overline{\Omega})$ have positive lower bounds, and $q\in L^\infty(\Omega)$ is nonnegative.
\end{ass}

For the spatially-dependent variable order $\alpha:\Omega\to(0,1)$ of time derivation $\partial_t^{\alpha(x)}U$ in problem \eqref{eq:2D3D:ibvp}, we assume the following condition, which is needed for the well-posedness of the direct problem (see e.g. \cite{HJK:2025:ISDVO,Ki1,Kian:2018:TFD} for more details).
\begin{ass}
\label{ass:alpha:basic}
$\alpha:\Om\to(0,1)$ is a measurable function satisfying $\overline{\alpha}:=\operatorname{esssup}\{\alpha(x)\,:\,x\in\Om\}\in(0,1)$, $\underline{\alpha}:=\operatorname{essinf}\{\alpha(x)\,:\,x\in\Om\}\in (0,1)$ and $\overline{\alpha}< 2\underline{\alpha}$.
\end{ass}

Throughout we assume the following condition on the boundary excitation $g$. This choice is mainly motivated by the forward problem \eqref{eq:2D3D:ibvp}. Similar conditions were imposed in related works \cite{IkehataKian:2023,HJK:2025:ISDVO,Kian:2018:TFD}. The extra regularity condition $\varphi_k\in W^{2-\frac{1}{r},r}(\partial\Omega)$, $k=3,\dots,M$, will be used in Theorem \ref{theorem:analyticity}  for proving the  regularity  $U\in C([0,+\infty);C^1(\overline{\Omega}))$ of the solution $U$ of problem \eqref{eq:2D3D:ibvp}.
\begin{ass}\label{ass:g}
The Dirichlet boundary condition $g$ takes the form $g(t,x)=\sum_{k=3}^M t^{k} \varphi_k(x)$ for all $(t,x)\in(0,\infty)\times\partial\Omega$ for some $M\ge3$, $r\in(d,6]$ and $\varphi_k\in W^{2-\frac{1}{r},r}(\partial\Omega)$, $k=3,\dots,M$, with $\varphi_M\not\equiv0$.
\end{ass}

\subsection{Uniqueness for piecewise-constant orders}
First we consider a piecewise constant variable order $\alpha$, which satisfies the following condition.
\begin{ass}\label{ass:stratified:structure}
The variable order $\alpha$ is a piecewise constant function of the form
\begin{equation}\label{alph}\alpha=\sum_{j=0}^n\alpha_j\mathds{1}_{\Omega_j},\end{equation}
where $\{\alpha_j\}_{j=1}^n\subset(0,1)$ and the subdomains $\{\Omega_j\}_{j=0}^n$ are pairwise disjoint, ${\Omega}=\bigcup_{j=0}^n{\Omega_j}$, and each $\Omega_j$ is measurable and has a positive Lebesgue measure in $\mathbb{R}^d$.
\end{ass}

The sequence of subdomains $\{\Omega_j\}_{j=0}^n$ satisfying Assumption \ref{ass:stratified:structure} can always be rearranged and merged so that $\{\alpha_k\}_{k=0}^n$ is  strictly increasing. We shall assume this arrangement whenever we adopt Assumption \ref{ass:stratified:structure}. Following Remark \ref{remark:Uhat:C1lambda} of Section 3, for $g$ given by Assumptions \ref{ass:alpha:basic}, the solution $U$ of
\eqref{eq:2D3D:ibvp} is lying in $C([0,+\infty);C^1(\overline{\Omega}))$ and, for any $x_0\in\partial\Omega$, the map $t\mapsto \partial_\nu U(t,x_0)$ is well defined and lying in $C([0+\infty))$.
Using these properties, we can state the following uniqueness results under one-point boundary measurement. The proof of the theorem is given in Section \ref{sec:inverse}.

\begin{theorem}\label{theorem:flux:lowfreq:asympt:onept}For $i=1,2$, fix $\alpha^i$ satisfying
Assumptions \ref{ass:alpha:basic} and \ref{ass:stratified:structure} with $n=n^i$, $\alpha_j=\alpha_j^i$, $\Omega_j=\Omega_j^i$, $j=0,\ldots,n^i$.
Assume also that $\rho^i$, $\sigma^i$ and $q^i$ fulfill Assumption \ref{ass:coef} and $g^i$ satisfies Assumption \ref{ass:g} with $M=M^i$ and $\phi_j=\phi_j^i$, $j=3,\ldots,M^i$. Let $U^i$, $i=1,2$, be the solution of problem \eqref{eq:2D3D:ibvp} with $\alpha=\alpha^i$, $\rho=\rho^i$, $\sigma=\sigma^i$, $q=q^i$ and $g=g^i$ for $i=1,2$.
Suppose that $\varphi_{M^1}^1$ and $\varphi_{M^2}^2$ are not sign-changing and fix any $x_0\in\partial\Omega$. Then the following statements hold.
\begin{itemize}
\item[{\rm(a)}] If there exists some sequence $\{t_k\}_{k=1}^\infty$ of distinct numbers accumulating to a positive number such that \begin{equation}\label{t1a} \sigma^1(x_0)\p_\nu  U^1(t_k,x_0)=\sigma^2(x_0)\p_\nu  U^2(t_k,x_0),\quad \forall k\in \mathbb{N},
\end{equation}
then there holds
\begin{equation}\label{t1b} \{\alpha^1(x)\,:\,x\in\Om\}=\{\alpha^2(x)\,:\, x\in\Om\}.
\end{equation}
\item[{\rm(b)}] Additionally, suppose that $\rho^1=\rho^2$, $\sigma^1=\sigma^2$, $q^1=q^2$ and $g^1=g^2$ and let one of the following conditions be fulfilled:
 \begin{itemize}
\item[{\rm(i)}] $\alpha^1\ge\alpha^2$;
\item[{\rm(ii)}]  $\Omega_j^1\subset\Omega_j^2$ or $\Omega_j^2\subset\Omega_j^1$, for $j=0,\ldots,\min\{n^1,n^2\}$.
\end{itemize}
Then the condition \eqref{t1a} implies that $\alpha^1=\alpha^2$ almost everywhere in $\Om$.
\end{itemize}
\end{theorem}

To the best of our knowledge, in Theorem \ref{theorem:flux:lowfreq:asympt:onept} we obtain the first result of a unique multidimensional determination
of a fractional time derivative of variable order $\alpha(x)$ from a single point flux measurement. Indeed, all the existing results have been either stated with an infinite number of measurements (i.e., Dirichlet-to-Neumann map) \cite{Kian:2018:TFD} or one single measurement on the full boundary \cite{IkehataKian:2023}. The works \cite{JinKian:2021procA,Jin:2022:ROD,JinKian:2023cms} proved the uniqueness of recovering one or multiple constant orders in the time fractional derivative, as well as the lower and upper bounds of the distributed order, from the flux measurement at one point on the boundary in the multi-dimensional setting. Our results can be seen as an extension of recent contributions \cite{JinKian:2021procA,Jin:2022:ROD,JinKian:2023cms}. In contrast to these last works \cite{JinKian:2021procA,Jin:2022:ROD,JinKian:2023cms}, we investigate the recovery of a map which needs much more information. It is worth emphasizing that we utilize only the measurement at one point on the boundary $\partial\Omega$, which makes the inverse problem far more challenging than the one-dimensional case, for which $\partial\Omega$ consists of only two points (cf. \cite{HJK:2025:ISDVO}). Since the flux data is only available at one boundary point, the problem is highly challenging, and we impose the monotonicity condition to overcome the challenge. The monotonicity condition is likely necessary to ensure the unique recovery of the $d$-dimensional target $\alpha$ from the one-dimensional flux data at one boundary point $x_0\in\partial\Omega$.

Our analysis is based on asymptotic properties of the Laplace transform $\partial_\nu \widehat{U}(p,x)$ in time of the flux $\partial_\nu U(t,x)$ at the frequency $p=0$ in Theorem \ref{theorem:flux:lowfreq:asympt:general}. In the analysis, we use the theory of resolvent operators defined on $L^r(\Omega)$ spaces, with $r> d$,  combined with Sobolev embedding properties in order to extract information from one single point measurement. This argument extends and greatly simplifies the gluing theory of \cite{HJK:2025:ISDVO} in dimension one where the analysis was based on Fourier series representation and a delicate analysis of a certain parameterized triangular system. This new approach allows much greater generality, especially dealing with a more general variable order $\alpha$.

\subsection{Uniqueness and stability in the generalized case}
We will state novel uniqueness, semi-continuity and conditional stability results in the generalized case in which the variable order $\alpha$ is not necessarily piecewise constant. First, we prove that any measurable variable order $\alpha$ satisfying Assumption \ref{ass:alpha:basic} can be uniquely determined by the one-point boundary measurement.
\begin{theorem}
\label{theorem:uniqueness:measfuns}
For $i=1,2$, fix $\alpha^i$ satisfying
Assumptions \ref{ass:alpha:basic}. Assume also that $\rho^i$, $\sigma^i$ and $q^i$ fulfill Assumption \ref{ass:coef} and $g^i$ satisfies Assumption \ref{ass:g} with $M=M^i$ and and $\phi_j=\phi_j^i$, $j=3,\ldots,M^i$. Let $U^i$ be the solution to problem \eqref{eq:2D3D:ibvp} with $\alpha=\alpha^i$, $\rho=\rho^i$, $\sigma=\sigma^i$, $q=q^i$ and $g=g^i$ for $i=1,2$.
Suppose that $\varphi_{M^1}^1$ and $\varphi_{M^2}^2$ are not sign-changing.
Fix any $x_0\in\partial\Omega$.
If there exists some sequence $\{t_k\}_{k=1}^\infty$ of distinct numbers accumulating to a positive number such that $$\sigma^1(x_0)\p_\nu  U^1(t_k,x_0)=\sigma^2(x_0)\p_\nu  U^2(t_k,x_0),\quad \forall k\in \mathbb{N}$$
then the following two statements hold:
\begin{itemize}
\item[\rm(a)] $\underline{\alpha^1}:=\operatorname{essinf}\{\alpha^1(x)\,:\,x\in\Om\}=\underline{\alpha^2}:=\operatorname{essinf}\{\alpha^2(x)\,:\,x\in\Om\}.$

\item[\rm(b)] 	Additionally	suppose  that $\rho^1=\rho^2$, $\sigma^1=\sigma^2$, $q^1=q^2$ and $g^1=g^2$ hold. If $\alpha^1\ge\alpha^2$ almost everywhere in $\Om$, then $\alpha^1=\alpha^2$ almost everywhere in $\Om$.
\end{itemize}
\end{theorem}

Part (a) of Theorem \ref{theorem:uniqueness:measfuns} gives the unique determination of the lower bound $\underline{\alpha}$ of the variable order $\alpha(x)$ from the one-point flux measurement of the variable-order subdiffusion model \eqref{eq:2D3D:ibvp} in an unknown medium (unknown coefficients and boundary condition), whereas part (b) gives the uniqueness of the full recovery of a general spatially variable order $\alpha$ based on the one-point boundary flux observation, under the monotonicity condition $\alpha^1\geq \alpha^2$. To the best of our knowledge, these are the first results in that category for a variable order that might not be piecewise constant.

Let us observe that the monotonicity condition imposed in Theorem \ref{theorem:uniqueness:measfuns} might be compared to conditions imposed in the inverse potential problem for the Schr\"{o}dinger equation. For the unique identifiability of the potential in the Schr\"{o}dinger equation in two- or three-dimension, we refer interested readers to  the work \cite[Corollary 5.2]{Harrach:2019:MLU} for the uniqueness under infinitely many measurements (DtN map) and the work \cite[Theorem 4.1]{Miao:2023:nonuniqueness} for the nonuniqueness under one measurement. Kim and Yamamoto \cite{Kim:2004:uniqueness,Kim:2003:uniqueness} proved several results on the unique determination of the support of the potential in the Schr\"{o}dinger equation using one boundary measurement and extra geometric information. In sharp contrast, Theorem \ref{theorem:uniqueness:measfuns}(b) gives the unique identifiability of the full information of general measurable function under one point flux measurement, and thus, it is not covered by these existing results.

Next, we prove a semi-continuity property with respect to   one point flux data, which significantly strengthens  the unique identifiability in Theorem \ref{theorem:uniqueness:measfuns} (b). For this purpose, let us observe that in view of the estimate \eqref{lala} in Remark \ref{remark:Uhat:C1lambda} of Section 3, for any $x_0\in\partial\Omega$, the map $t\mapsto e^{-t/2}\partial_\nu U(t,x_0)$ belongs to $L^1(0,+\infty)$. Using this property, we can state a semi-continuity property as follows.

\begin{theorem}
\label{theorem:stability}
Let Assumption \ref{ass:coef} hold, and let $g$ satisfy Assumption \ref{ass:g} and its Laplace transform $\widehat{g}(1,x)\equiv\sum_{k=3}^M k!\varphi_k(x)$ be not everywhere-vanishing nor sign-changing on the boundary $\p\Omega$.
Fix $\alpha^1$ and the sequence of functions $\{\alpha^{2,\ell}:\ \ell\in\mathbb N\}$ satisfying Assumption
\ref{ass:alpha:basic} with $\alpha=\alpha^1$ or $\alpha=\alpha^{2,\ell}$, $\ell\in\mathbb N$. Assume also that $\alpha^{2,\ell}\geq \alpha^1$, $\ell\in\mathbb N$, and the sequence $\{\alpha^{2,\ell}:\ \ell\in\mathbb N\}$ is $L^r(\Omega)$-convergent. For all $\ell\in\mathbb N$, fix $U^1$ {\rm(}respectively $U^{2,\ell}${\rm)} the solution to problem \eqref{eq:2D3D:ibvp} with $\alpha=\alpha^1$ {\rm(}respectively $\alpha=\alpha^{2,\ell}${\rm)}.
Fix any $x_0\in\partial\Omega$.
Then the condition
$$\lim_{\ell\to\infty}\int_0^\infty|\p_\nu  U^1(t,x_0)-\p_\nu  U^{2,\ell}(t,x_0)| (te^{-t})\d t=0$$
implies that
$$\lim_{\ell\to\infty}\|\alpha^1-\alpha^{2,\ell}\|_{L^r(\Omega)}= 0.$$
\end{theorem}

\begin{remark}
One can remove the assumption that $\alpha^{2,\ell}$ is convergent in $L^r(\Om)$ if we restrict the variable order $\alpha$ to be in a compact set of measurable functions.
Also, choosing $\alpha^{2,\ell}$ to be the same for all $\ell\in\mathbb{N}$ in Theorem \ref{theorem:stability} gives another version of the uniqueness result that is similar to Theorem \ref{theorem:uniqueness:measfuns} {\rm(}b{\rm)} but with a different assumption on the Dirichlet condition $g$.
\end{remark}

The continuity property in Theorem \ref{theorem:stability} involves the weight $te^{-t}$, which arises naturally from the derivative of the Laplace transform $\partial_\nu \widehat{U}(p,x)$ with respect to the frequency $p$, cf. \eqref{eq:limit:diffp1} in Lemma \ref{lem:int-repres}.
Theorem \ref{theorem:stability} gives a continuity property, in addition to the uniqueness result of
Theorem \ref{theorem:uniqueness:measfuns}(b). Such continuity property is weaker than the classical conditional stability, but it can still quantify the uniqueness result. The proof of the continuity result is still based on suitable asymptotic properties of the Laplace transform in time  of the solution of problem \eqref{eq:2D3D:ibvp}. However, the asymptotic properties are not considered at the frequency $p=0$ but at $p=1$. This new approach allows strengthening our uniqueness results into continuity results. Note that the monotonicity condition enables the use of the strong maximum principle and Hopf's lemma, and plays a crucial role in the analysis.

In order to upgrade the semi-continuity result of Theorem \ref{theorem:stability} to a conditional stability result, we need extra boundary data. Namely, given the full boundary data, we obtain the following conditional Lipschitz stability under the monotonicity assumption on the variable order $\alpha$.
We again use the estimate \eqref{lala}, which appears in Remark \ref{remark:Uhat:C1lambda} of Section 3, to prove that the map $(t,x)\mapsto e^{-t/2}\partial_\nu U(t,x)$ belongs to $L^1(0,+\infty; L^2(\partial\Omega))$.

\begin{theorem}
\label{theorem:stability:fullbdydata}
Let Assumptions \ref{ass:coef} and \ref{ass:g} be fulfilled and fix  $\alpha^i$, $i=1,2$, satisfying Assumption \ref{ass:alpha:basic} for $\alpha=\alpha^i$ and the monotonicity condition $\alpha^1\ge\alpha^2$ almost everywhere in the domain $\Omega$.
Let $U^i$, $i=1,2$, be the solution to problem \eqref{eq:2D3D:ibvp} with $\alpha=\alpha^i$and suppose that there exists a constant $g_0>0$ such that $|\widehat{g}(1,x)|\equiv|\sum_{k=3}^M k!\varphi_k(x)|\ge g_0$ for all $x\in\partial\Omega$.
Then, there exists $C>0$,  depending only on $\sigma$, $\rho$, $\Omega$, $q$ and $g$, such that
\begin{equation}\label{esti}
    \|\alpha^1-\alpha^{2}\|_{L^1(\Omega)}\le C\int_{\p\Om}\left|\int_{0}^\infty \left(\p_\nu  U^1(t,x)-\p_\nu  U^{2}(t,x)\right)(te^{-t})\d t\right|\d s.\end{equation}
\end{theorem}

To the best of our knowledge, this is the first stability result for recovering the variable order $\alpha$. In the case of one constant order, the only work we are aware of is due to Li et al \cite{LiHuangYamamoto:2020}, who proved the Lipschitz stability of recovering the constant order $\alpha$ from the measurement of the solution $U(t_0,x_0)$ at one interior point $x_0\in\Omega$ at a fixed time instance $t_0$.

Note that the conditional stability result of Theorem \ref{theorem:stability:fullbdydata} is only subject to the monotonicity condition imposed on the unknown parameters $\alpha^1$ and $\alpha^2$. However, no extra a priori upper or lower bounds are imposed to these unknown parameters.

\subsection{Illustration of the theorems}
\begin{figure}[hbt!]
\centering\setlength{\tabcolsep}{0pt}
\begin{tabular}{ccccc}
	\includegraphics[width=0.23\linewidth]{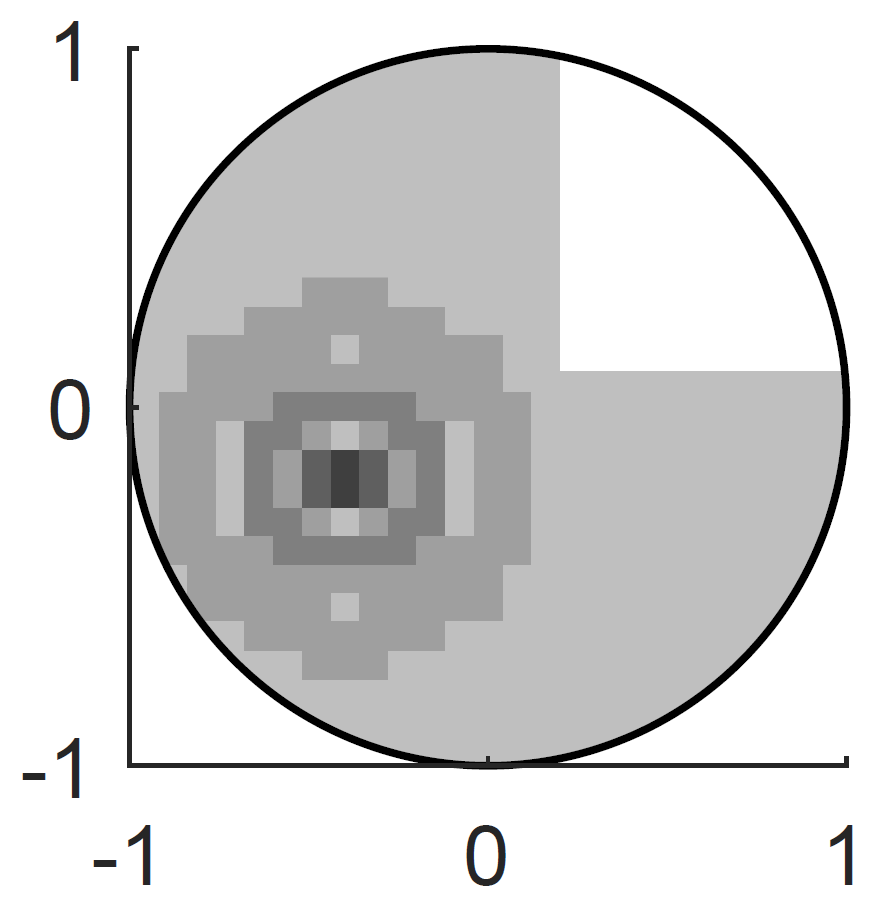}&
\includegraphics[width=0.23\linewidth]{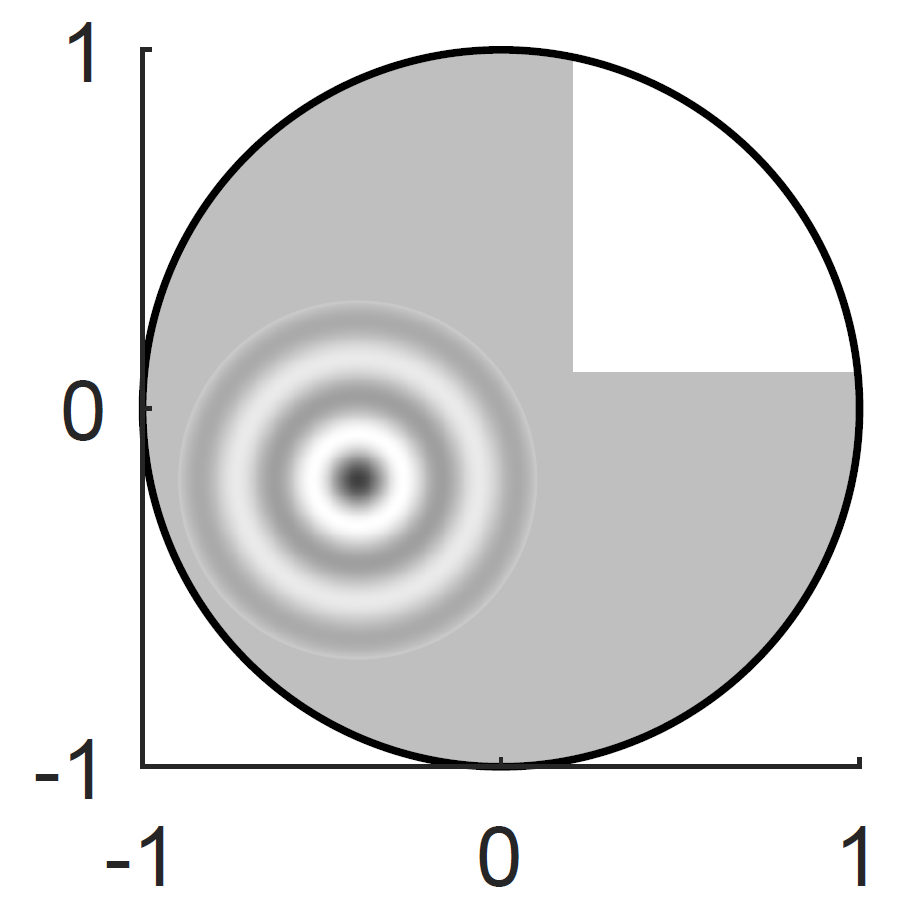}&
\includegraphics[width=0.23\linewidth]{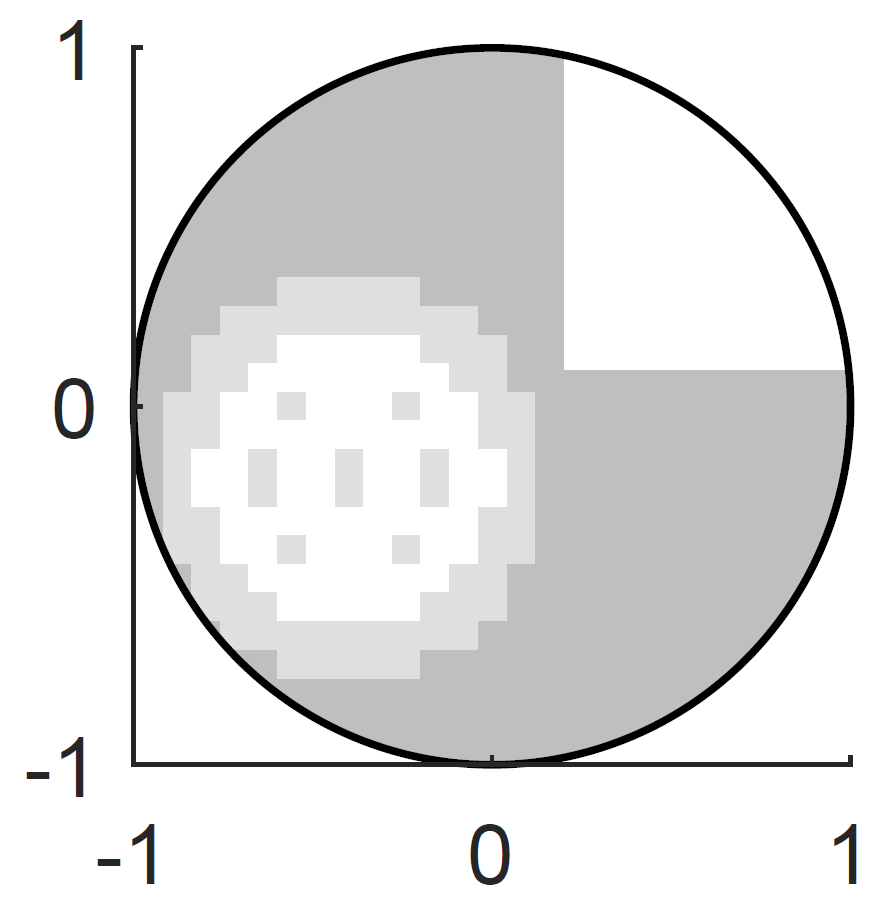}&
\includegraphics[width=0.23\linewidth]{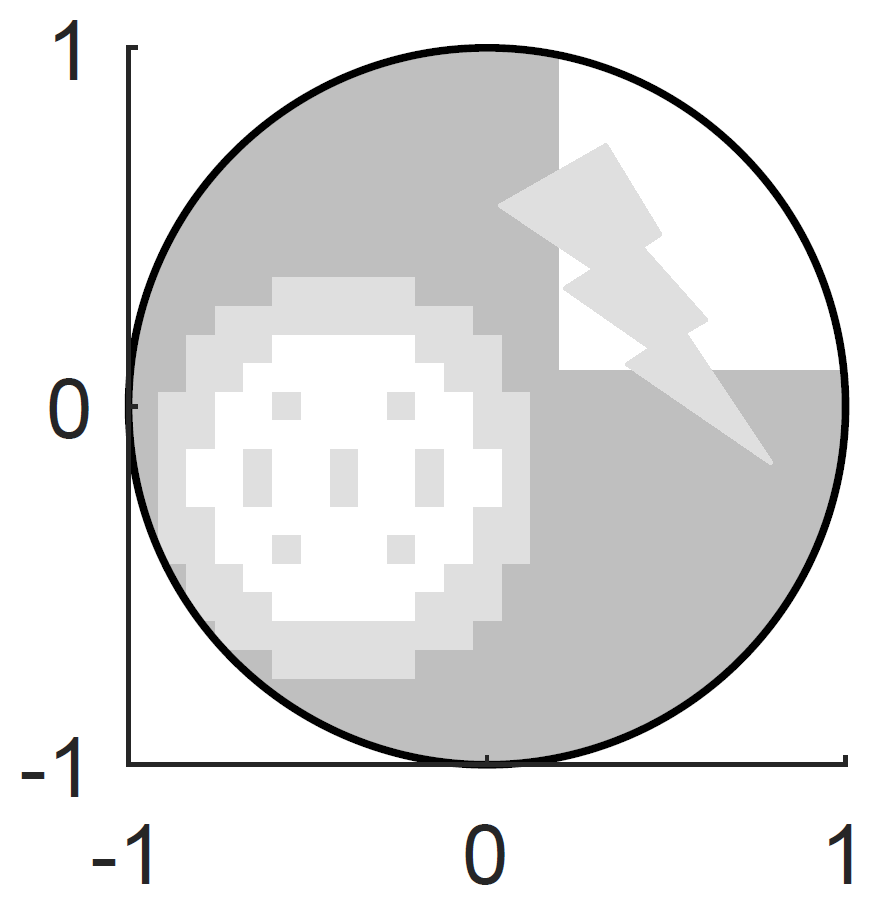}&
\includegraphics[width=0.062\linewidth]{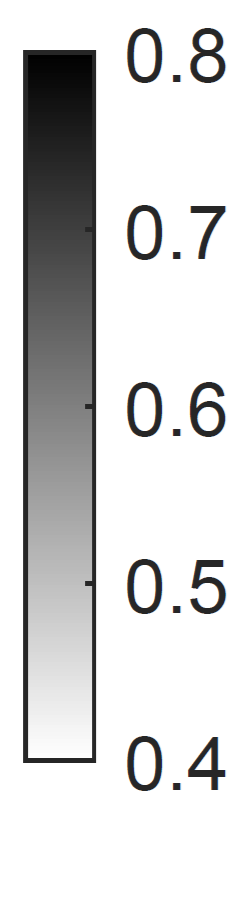}  \\
(a) $\alpha=\alpha^1$ & (b) $\alpha=\alpha^2$ & (c) $\alpha=\alpha^3$ & (d) $\alpha=\alpha^4$ &
\end{tabular}
\caption{\label{fig:all:examples}
The contour plots of four variable orders $\alpha:B_1(0)\to[0.4,0.8]$ satisfying Assumption \ref{ass:alpha:basic}. The functions $\alpha^1$, $\alpha^3$ and $\alpha^4$ are piecewise constant functions satisfying Assumption \ref{ass:stratified:structure}, and there hold the monotonicity relation $\alpha^1\ge \alpha^2\ge \alpha^3$ almost everywhere.
$\alpha^3$ and $\alpha^4$ satisfy the second assumption of part (b) of Theorem \ref{theorem:flux:lowfreq:asympt:onept}.}
\end{figure}

Now we briefly illustrate various situations covered by the theorems. Consider the unit disk $\Omega = B_1(0)$ centered at the origin.
Fig. \ref{fig:all:examples} shows four variable orders $\{\alpha^i\}_{i=1}^4$, where $\alpha^1$, $\alpha^3$ and $\alpha^4$ are piecewise constant, and $\alpha^2$ belongs to the generalized case.
By the construction, we have $\alpha^1\ge\alpha^2\ge\alpha^3$ almost everywhere. This monotonicity relation and Theorem \ref{theorem:uniqueness:measfuns} (b) imply that one point flux measurement can distinguish $\alpha^1$, $\alpha^2$ and $\alpha^3$.
However, since there is no monotonicity relation between $\alpha^i$ and $\alpha^4$ for any $i\in\{1,2,3\}$, so Theorem \ref{theorem:uniqueness:measfuns} cannot be used to distinguish them.
Instead, $\alpha^4$ can be distinguished from $\alpha^1$ and $\alpha^3$ by one point flux measurement in view of Theorem \ref{theorem:flux:lowfreq:asympt:onept} as follows. For $i\in\{3,4\}$, since $0.7$ is in the range of $\alpha^1$ but not in the range of $\alpha^i$, $\alpha^1$ and $\alpha^i$ are distinguished by one point flux measurement via part (a) of Theorem \ref{theorem:flux:lowfreq:asympt:onept}(even for unknown medium properties $\rho$, $\sigma$ and $q$ and unknown excitation $g$).
Part (b) of Theorem \ref{theorem:flux:lowfreq:asympt:onept} and the second assumption on $\alpha$ indicate that one point flux measurement is enough to distinguish $\alpha^3$ and $\alpha^4$.

\section{Preliminaries on the direct problem}\label{sec:direct}
In this section, we develop several key analytic properties of the direct problem \eqref{eq:2D3D:ibvp}, which play a crucial role in the analysis of the inverse problem in Section \ref{sec:inverse}.

\subsection{Well-posedness of the direct problem}
First we discuss the well-posedness of problem \eqref{eq:2D3D:ibvp} using $L^r(\Omega)$ resolvent estimates.
Let the coefficients $\rho$, $\sigma$ and $q$ satisfy Assumption \ref{ass:coef}. Then, for $\kappa\in[2,+\infty)$, we define an elliptic operator $\mathcal{A}$ on the space $W^{2,\kappa}(\Omega)$ by
$$\mathcal{A}[v](x):=\rho(x)^{-1}(-\nabla\cdot(\sigma(x)\nabla v(x))+q(x)v(x)),\quad\forall  v\in W^{2,\kappa}(\Omega).$$

The next result gives basic properties of the elliptic operator $\mathcal{A}$.
\begin{lemma}\label{lemma:iso:Apluspalpha}Let Assumptions \ref{ass:coef} and \ref{ass:alpha:basic} hold. For each fixed $p\geq 0$ and $\kappa\in[2,+\infty)$, the map $u\mapsto((\mathcal{A}+p^{\alpha})u,u|_{\partial\Omega})$ defines an isomorphism from $W^{2,\kappa}(\Om)$ onto $L^{\kappa}(\Om)\times W^{2-1/\kappa,\kappa}(\partial\Omega)$. Moreover, for all $\kappa\in [2,6]$, there exists a constant $C>0$, depending only on $\mathcal A$, $\rho$, $\kappa$ and $\Omega$, such that for all $p\geq0$, there holds
\begin{equation}\label{l1a}\|v\|_{L^\kappa(\Omega)}\leq C\|(\mathcal{A}+p^{\alpha})v\|_{L^\kappa(\Omega)},\quad \forall v\in W^{2,\kappa}(\Omega)\cap H^1_0(\Omega).
\end{equation}
\end{lemma}
\begin{proof}
Without loss of generality, we consider only real valued functions.
By \cite[Theorem 2.4.2.5 and Remark 2.5.1.2]{Grisvard:1985:EPND}, the map $u\mapsto(\rho(\mathcal{A}+p^{\alpha})u,u|_{\partial\Omega})$ defines an isomorphism from $W^{2,\kappa}(\Om)$ onto $L^{\kappa}(\Om)\times W^{2-1/\kappa,\kappa}(\partial\Omega)$.
Since $\rho$ is bounded from below and above by positive numbers, the same assertion holds for the map $u\mapsto((\mathcal{A}+p^{\alpha})u,u|_{\partial\Omega})$. To derive the estimate \eqref{l1a}, by fixing $v\in W^{2,\kappa}(\Omega)\cap H^1_0(\Omega)$, $p\geq0$, and using Poincar\'e inequality, we deduce
$$\begin{aligned}
&\|v\|_{H^1(\Omega)}^2\leq C\int_\Omega \sigma |\nabla v|^2{\rm d}x=C\int_\Omega v(-\nabla\cdot\sigma \nabla v){\rm d}x\\
\leq &C\int_\Omega \rho v(\mathcal{A}v+p^{\alpha}v){\rm d}x\leq C \|v\|_{L^\kappa(\Omega)}\|\mathcal{A}v+p^{\alpha}v\|_{L^\kappa(\Omega)},\end{aligned}$$
where $C>0$ depends only on $\mathcal A$, $\rho$, $\kappa$ and $\Omega$. Combining this estimate with the Sobolev embedding $H^1(\Omega)\hookrightarrow L^\kappa(\Omega)$ for $\kappa\in[2,6]$ yields the desired estimate \eqref{l1a}.
\end{proof}

Below, for all $\kappa\in[2,+\infty)$, we denote by $\mathbb A_\kappa$ the unbounded operator acting on $L^\kappa(\Omega)$ with its domain $D(\mathbb A_\kappa)=H_0^1(\Om)\cap W^{2,\kappa}(\Om)$ defined by
$$\forall v\in D(\mathbb A_\kappa),\quad \mathbb A_\kappa v=\mathcal A v.$$
For two Banach spaces ${X}$ and ${Y}$, $\mathcal{B}(X)$ denotes the operator norm from $X$ to itself, and $\mathcal{B}(X,Y)$ the operator norm from $X$ to $Y$. Then the following resolvent estimate holds.
\begin{lemma}\label{lem:Lr-resolvent}
Let Assumptions \ref{ass:coef} and \ref{ass:alpha:basic} hold. Fix $r\in (d,6]$. Then for all $\psi\in(0,\pi)$, the following $L^r(\Omega)$ resolvent estimate holds
\begin{equation*}
   \|(\mathbb{A}_r+(se^{i\varphi})^{\alpha(x)})^{-1}\|_{\mathcal{B}(L^r(\Omega),W^{2,r}(\Omega))}\leq C\max(s^{\underline{\alpha}-2\overline{\alpha}},s^{3\overline{\alpha}-2\underline{\alpha}}), \quad s>0, \varphi \in (-\psi,\psi),
\end{equation*}
where the constant $C>0$ is independent of $s>0$.
\end{lemma}
\begin{proof}
By \cite[Proposition 2.1]{Kian:2018:TFD} (see also \cite[formula (2.8)]{Ki1}), we can prove that, for all $p\in\mathbb C\setminus\mathbb R_-$, $\mathbb A_2+p^{\alpha(x)}$ is boundedly invertible in $L^2(\Omega)$, $(\mathbb A_2+p^{\alpha(x)})^{-1}$ maps $L^2(\Omega)$ to $H^2(\Omega)\cap H^1_0(\Omega)$ and, moreover, for all $\psi\in(0,\pi)$, there holds
\begin{align*}
  &\|(\mathbb A_2+(se^{i\phi})^{\alpha(x)})^{-1}\|_{\mathcal B(L^2(\Omega),H^2(\Omega))}\\
  \leq &C(\|\mathbb A_2(\mathbb A_2+(se^{i\phi})^{\alpha(x)})^{-1}\|_{\mathcal B(L^2(\Omega))}+\|(\mathbb A_2+(se^{i\phi})^{\alpha(x)})^{-1}\|_{\mathcal B(L^2(\Omega))})\\
  \leq &C(\|\mathrm{id}-(se^{i\phi})^{\alpha(x)}(\mathbb A_2+(se^{i\phi})^{\alpha(x)})^{-1}\|_{\mathcal B(L^2(\Omega))}+\|(\mathbb A_2+(se^{i\phi})^{\alpha(x)})^{-1}\|_{\mathcal B(L^2(\Omega))})\\
  \leq & C(1+(s^{\overline{\alpha}}+1)\|(\mathbb A_2+(se^{i\phi})^{\alpha(x)})^{-1}\|_{\mathcal B(L^2(\Omega))}) \leq  C(1+(s^{\overline{\alpha}}+1)\max(s^{\underline{\alpha}-2\overline{\alpha}},s^{\overline{\alpha}-2\underline{\alpha}}))\\
  \leq &C\max(s^{\underline{\alpha}-2\overline{\alpha}},s^{2(\overline{\alpha}-\underline{\alpha})}),\quad s>0,\ \phi\in(-\psi,\psi),
\end{align*}
where $C>0$ is independent of $s>0$, and $\mathrm{id}$ denotes the identity operator. Using the Sobolev embedding $H^2(\Omega)\hookrightarrow L^r(\Omega)$, we deduce that
for all $p\in\mathbb C\setminus\mathbb R_-$, $\mathbb A_r+p^{\alpha(x)}$ is boundedly invertible and, for all $\psi\in(0,\pi)$, we have
\begin{align*}
\|(\mathbb A_r+(se^{i\phi})^{\alpha(x)})^{-1}\|_{\mathcal B(L^r(\Omega))}&\leq C \|(\mathbb A_2+(se^{i\phi})^{\alpha(x)})^{-1}\|_{\mathcal B(L^2(\Omega),H^2(\Omega))}\\
&\leq C\max(s^{\underline{\alpha}-2\overline{\alpha}},s^{2(\overline{\alpha}-\underline{\alpha})}),\quad s>0,\ \phi\in(-\psi,\psi),
\end{align*}
where $C>0$  is independent of $s>0$. By combining this estimate with Lemma \ref{lemma:iso:Apluspalpha}, we obtain
\begin{align*}
&\|(\mathbb A_r+(se^{i\phi})^{\alpha(x)})^{-1}\|_{\mathcal B(L^r(\Omega),W^{2,r}(\Omega))}\\
\leq& C(\|\mathbb A_r(\mathbb A_r+(se^{i\phi})^{\alpha(x)})^{-1}\|_{\mathcal B(L^r(\Omega))}+\|(\mathbb A_r+(se^{i\phi})^{\alpha(x)})^{-1}\|_{\mathcal B(L^r(\Omega))})\\
\leq &C(\|\mathrm{id}-(se^{i\phi})^\alpha(\mathbb A_r+(se^{i\phi})^{\alpha(x)})^{-1}\|_{\mathcal B(L^r(\Omega))}+\|(\mathbb A_r+(se^{i\phi})^{\alpha(x)})^{-1}\|_{\mathcal B(L^r(\Omega))})
\\
\leq &C\max(s^{\underline{\alpha}-2\overline{\alpha}},s^{3\overline{\alpha}-2\underline{\alpha}}),\quad s>0,\ \phi\in(-\psi,\psi),
\end{align*}
with $C>0$  independent of $s>0$. This completes the proof of the lemma.
\end{proof}

Using these preliminary results, we can now prove the following extension of  \cite[Proposition 3.1]{Kian:2018:TFD} and \cite[Theorem 3.1]{HJK:2025:ISDVO} about the regularity and time analyticity of the solution $U(t,x)$ of  \eqref{eq:2D3D:ibvp}.
\begin{theorem}
\label{theorem:analyticity}
Under Assumptions \ref{ass:coef}, \ref{ass:alpha:basic} and \ref{ass:g}, problem \eqref{eq:2D3D:ibvp} admits a unique solution\\ $U\in C([0,+\infty);W^{2,r}(\Omega))\cap C^1([0,+\infty);L^2(\Omega))$. Moreover, the solution $U(t,\cdot)$ is analytic in time $t\in(0,+\infty)$ as a map taking values in $W^{2,r}(\Omega)$, and its Laplace transform $\widehat{U}(p,\cdot):=\int_0^\infty U(t.x)e^{-pt}\d t$ is in $W^{2,r}(\Omega)$ for all $p>0$.
\end{theorem}
\begin{proof} Without loss of generality, we assume that Assumption \ref{ass:g} is fulfilled with $M=3$.
 In view of \cite[Theorem 2.4.2.5]{Grisvard:1985:EPND}, let $\Phi\in W^{2,r}(\Omega)$ be the unique solution of
\begin{equation}\label{BVP}
\left\{\begin{aligned}
- \nabla \cdot( \sigma \nabla\Phi) + q\Phi&=0,\quad\mbox{in } \Omega,\\
\Phi&=\phi_3,\quad\mbox{on }\partial\Omega.
\end{aligned}\right.
\end{equation}
Also let $\widehat W$ be defined by
$$\widehat W(p,\cdot)=-6p^{-4}(\mathbb A_r+p^{\alpha(x)})^{-1}p^{\alpha(x)}\Phi,\quad p\in\mathbb C\setminus\mathbb R_-.$$
Fix
$\theta\in(\frac{\pi}{2},\pi)$, $\delta>0$ and define a contour
$\gamma(\delta,\theta):=\gamma_-(\delta,\theta)\cup\gamma_0(\delta,\theta)\cup\gamma_+(\delta,\theta)\subset\mathbb{C}$, oriented  counterclockwise with
$\gamma_0(\delta,\theta):=\{\delta\, e^{i\beta}:\ \beta\in[-\theta,\theta]\}$, $\gamma_\pm(\delta,\theta)
:=\{s\,e^{\pm i\theta}: s\in[\delta,\infty)\}$.
By combining the argument of \cite[Proposition 3.1]{Kian:2018:TFD} with \cite[p. 16]{Ki1} (see also \cite[Theorem 3.1]{HJK:2025:ISDVO} for a similar argument), we can deduce that
problem \eqref{eq:2D3D:ibvp} admits a unique solution $U\in C([0,+\infty);H^2(\Omega))\cap C^1([0,+\infty);L^2(\Omega))$ defined by
$$U(t,x)=t^3\Phi(x)+\frac{1}{2i\pi}\int_{\gamma(\delta,\theta)}e^{t p} \widehat W(p,\cdot)\,{\rm d} p,\quad (t,x)\in[0,+\infty)\times\Omega.$$
By combining this representation with Lemma \ref{lem:Lr-resolvent}, in a similar way to \cite[Proposition 3.1]{Kian:2018:TFD} and \cite[Theorem 3.1]{HJK:2025:ISDVO}, we deduce that
$U\in C([0,+\infty);W^{2,r}(\Omega))$. The time analyticity can be deduced by combining the preceding argumentation with \cite[Lemma 3.2]{Kian:2018:TFD}.
\end{proof}

\begin{remark}[Regularity of $U$ and $\widehat{U}$]
\label{remark:Uhat:C1lambda} Recall that under Assumption \ref{ass:g}, we have $r\in (d,6]$. Combining  the Sobolev embedding theorem  \cite[Part II of Chapter 4.12]{Adams:2003:SS} with Theorem \ref{theorem:analyticity} {\rm(}under Assumptions \ref{ass:coef}, \ref{ass:alpha:basic} and \ref{ass:g}{\rm)} shows that the map $t\mapsto U(t,\cdot)$ is analytic with respect to $t>0$ as a map taking values in $W^{2,r}(\Omega)\hookrightarrow C^{1,r_0}(\overline{\Omega}),$ for all $0<r_0\le 1-\frac{d}{r}$. Therefore, for all $x_0\in\partial\Omega$, the map $t\mapsto \partial_\nu U(t,x_0)$ is analytic with respect to $t>0$ and continuous with respect to $t\geq0$. Moreover,
the Laplace transform $\widehat{U}(p,\cdot)$ of $U(t,\cdot)$ satisfies, for $p>0$,
\begin{equation}\label{eq:afterLaptrans}
\left\{\begin{aligned}
- \nabla \cdot( \sigma \nabla\widehat{U}(p,\cdot)) + (q + \rho p^{\alpha(x)})\widehat{U}(p,\cdot)&=0,\quad\mbox{in } \Omega,\\
\widehat{U}(p,\cdot)&=\widehat{g}(p,\cdot),\quad\mbox{on }\partial\Omega
\end{aligned}\right.
\end{equation}
and $\widehat{U}(p,\cdot)\in W^{2,r}(\Omega)\hookrightarrow C^{1,r_0}(\overline{\Omega}),$ for all $0<r_0\le 1-\frac{d}{r}$. Finally, by combining Lemma \ref{lem:Lr-resolvent} with  the argumentation of \cite[Proposition 3.1]{Kian:2018:TFD} and \cite[Lemma 2.1]{Ki1}, we can show that, for all $\epsilon>0$, there holds
\begin{equation}\label{lala}\int_0^{+\infty}e^{-\epsilon t}\|U(t,\cdot)\|_{C^1(\overline{\Omega})}{\rm d}t\leq C\int_0^{+\infty}e^{-\epsilon t}\|U(t,\cdot)\|_{W^{2,r}(\Omega)}{\rm d}t<\infty.
\end{equation}
\end{remark}

\subsection{Representations of the solution and flux in the Laplace domain}
In the rest of this work, we assume that $r\in(d,6]$, as specified in Assumption \ref{ass:g}. By Lemma \ref{lem:Lr-resolvent}, we deduce
\begin{equation}\label{estt1}\|(\mathbb A_r+p^{\alpha})^{-1}h\|_{W^{2,r}(\Om)}\le C\|h\|_{L^r(\Om)},\quad\forall h\in L^r(\Om),\end{equation}
where the constant $C>0$ is independent of $p\geq0$.

From Lemma \ref{lemma:iso:Apluspalpha}, we deduce that, for all $p\geq0$, $(\mathbb A_r+p^{\alpha})^{-1}$ is well-defined as a  bounded linear map from $L^{r}(\Om)$ to $H_0^1(\Om)\cap W^{2,r}(\Om)$. By Lemma \ref{lemma:iso:Apluspalpha} with $p=0$, we can also define the solution operator $\mathcal{S}:W^{2-1/r,r}(\partial\Omega)\to W^{2,r}(\Om)$ by $v_0\mapsto v$ for the unique solution $v$ of
\begin{equation}\label{eq:S:bvp:wholedomain}
\left\{
\begin{aligned}
-\nabla\cdot(\sigma\nabla v) + q v&=0,\quad \mbox{in }\Om,\\
v&=v_0,\quad\mbox{on }\partial\Omega.
\end{aligned}
\right.
\end{equation}

Now we can give an important solution representation. Let $\widehat{U}(p,x)$ be the Laplace transform of $U(t,x)$ in time $t$.
By Remark \ref{remark:Uhat:C1lambda}, we have $\widehat{U}(p,\cdot)\in W^{2,r}(\Omega)\hookrightarrow C^{1,r_0}(\overline{\Omega})$ for all $p>0$.
The next lemma can be alternatively stated as an eigenfunction expansion; see \cite[Lemma 3.2]{HJK:2025:ISDVO} for a piecewise constant $\alpha$ in the one-dimension case.
\begin{lemma}
\label{lemma:representation:resolvents}
Fix any $p>0$. Under Assumptions \ref{ass:coef}, \ref{ass:alpha:basic} and \ref{ass:g}, we have
\begin{align}
\widehat{U}(p,\cdot)&=\left({\rm id}-(\mathbb A_r+p^{\alpha})^{-1}p^{\alpha}\right)\mathcal{S}[\widehat{g}(p,\cdot)],\quad\mbox{in }\Om,
\label{eq:Uhat:repr}
\end{align}
where $(\mathbb A_r+p^{\alpha})^{-1}$ is a bounded linear operator from $L^{r}(\Omega)$ to $H_0^1(\Omega)\cap W^{2,r}(\Omega)$.
Moreover, for all $0<r_0\le 1-\frac{d}{r}$, there holds
\begin{equation}\label{eq:normalder:Uhat}
 \p_\nu \widehat{U}(p,\cdot)\big|_{\partial\Omega}=
\p_\nu \left({\rm id}-(\mathbb A_r+p^{\alpha})^{-1}p^{\alpha}\right)\mathcal{S}[\widehat{g}(p,\cdot)]\big|_{\partial\Omega}\in C^{0,r_0}(\partial\Omega).
\end{equation}
\end{lemma}
\begin{proof}
From the relations $\mathcal{A}\mathcal{S}\widehat{g}(p,\cdot)=0$ and $(\mathcal{A}+p^{\alpha})\widehat{U}(p,\cdot)=0$ in $\Om$, we obtain
$$(\mathcal{A}+p^{\alpha})\left[\widehat{U}(p,\cdot) - \mathcal{S}\widehat{g}(p,\cdot)\right] = -p^{\alpha}\mathcal{S}[\widehat{g}(p,\cdot)],\quad\mbox{in }\Omega.$$
This identity, the relation $\widehat{U}(p,\cdot) - \mathcal{S}\widehat{g}(p,\cdot)\in H_0^1(\Omega)$ and Lemma \ref{lemma:iso:Apluspalpha} imply the desired identity \eqref{eq:Uhat:repr}.
The Sobolev embedding $W^{2,r}(\Om)\hookrightarrow C^{1,r_0}(\overline{\Omega})$ \cite[Part II of Chapter 4.12]{Adams:2003:SS} shows the H\"{o}lder continuity of the normal derivative $\p_\nu \widehat{U}(p,\cdot)\big|_{\partial\Omega}$.
\end{proof}

\subsection{Asymptotic expansions}
Now we present two asymptotic expansions of the normal derivative $\p_\nu \widehat{U}(p,\cdot)$  using the expression \eqref{eq:normalder:Uhat}, one as $p\to 0^+$, and the other near $p=1$. These two expansions provide crucial technical tools for establishing the uniqueness and stability results stated in Section \ref{sec:main}.
The asymptotic expansion is based on the following Neumann series expansion.

\begin{lemma}\label{lemma:neumannseries}
Let Assumptions \ref{ass:coef} and \ref{ass:alpha:basic} hold. Fix any $N\in\mathbb{N}$. There exist $C>0$ and $p_0\in(0,1)$ such that for all $p\in(0,p_0)$ and $f\in L^r(\Om)$:
\begin{equation}\label{eq:bound:Neumannseries}
\left\|(\mathbb A_r+p^{\alpha})^{-1}f-\mathbb A_r^{-1}\sum_{i=0}^{N-1}(-p^{\alpha}\mathbb A_r^{-1})^i f\right\|_{W^{2,r}(\Om)}\le Cp^{N\underline{\alpha}}\|f\|_{L^r(\Om)}.
\end{equation}
\end{lemma}
\begin{proof}
Recall that $(\mathbb A_r + p^\alpha)^{-1}$ is well-defined as a bounded linear map from $L^{r}(\Om)$ to $H_0^1(\Om)\cap W^{2,r}(\Om)$ for all $p\ge0$. Since  by Assumption \ref{ass:coef}, $q$ is nonnegative, $0$ is not in the spectrum of the operator $\mathbb A_r$. Thus we can fix $C_A:=\|\mathbb A_r^{-1}\|_{L^r(\Om)\to W^{2,r}(\Om)}<\infty$.
Then, we have
$$\|p^{\alpha} \mathbb A_r^{-1}\|_{L^r(\Om)\to L^r(\Om)}\le C_A\|p^{\alpha}\|_{L^\infty(\Om)}=C_A p^{\underline{\alpha}},\quad \forall p\in(0,1).$$
Fix any $p_0\in(0,1)$ satisfying $C_A p_0^{\underline{\alpha}}<1/2$.
Then, for all $p\in(0,p_0)$, the operator ${\rm id}+p^{\alpha}\mathbb A_r^{-1}$ is invertible on $L^r(\Om)$ and satisfies the following estimate for all $p\in(0,p_0)$:
\begin{align*}
&\quad\left\|({\rm id}+p^{\alpha}\mathbb A_r^{-1})^{-1}-\sum_{i=0}^{N-1}(-p^\alpha \mathbb A_r^{-1})^i\right\|_{L^r(\Om)\to L^r(\Om)}\le \sum_{i\ge N}\|p^\alpha \mathbb A_r^{-1}\|_{L^r(\Om)\to L^r(\Om)}^i\\
&
\le \sum_{i\ge N} (C_A p^{\underline{\alpha}})^i = \frac{(C_A p^{\underline{\alpha}})^N}{1-C_A p^{\underline{\alpha}}}\le 2 C_A^N p^{N\underline{\alpha}}.
\end{align*}
Also, there holds the identity $(\mathbb A_r+p^{\alpha})^{-1}=\mathbb A_r^{-1}({\rm id}+p^{\alpha}\mathbb A_r^{-1})^{-1}$ for all $p\in(0,p_0)$.
Therefore, the desired estimate \eqref{eq:bound:Neumannseries} holds with $C=2C_A^{N+1}$.
\end{proof}

Note that similar to the case $p=0$, for the case $p=1$, the value of $p^\alpha$ is always $1$, which is independent of $\alpha$. This observation motivates the following Neumann series expansion near $p=1$.

\begin{lemma}\label{lemma:neumannseries:near:1}
Let Assumptions \ref{ass:coef} and \ref{ass:alpha:basic} hold. Fix any $N\in\mathbb{N}$. There exist  $C,\varepsilon>0$ satisfying the following estimate for all $p\in(1-\varepsilon,1+\varepsilon)$ and $f\in L^r(\Om)$:
\begin{equation}\label{eq:bound:Neumannseries:near:1}
\begin{aligned}
&\left\|(\mathbb A_r+p^{\alpha})^{-1}f-(\mathbb A_r+{\rm id})^{-1}\sum_{i=0}^{N-1}(({\rm id}-p^{\alpha})(\mathbb A_r+{\rm id})^{-1})^i f\right\|_{W^{2,r}(\Om)}\\
\le& C\|{\rm id}-p^\alpha\|_{L^\infty(\Om)}^{N}\|f\|_{L^r(\Om)}.
\end{aligned}
\end{equation}
\end{lemma}
\begin{proof}
The proof is similar to that of Lemma \ref{lemma:neumannseries}. We use the identity $$(\mathbb A_r+p^{\alpha})^{-1}=(\mathbb A_r+{\rm id})^{-1}({\rm id}+(p^{\alpha}-{\rm id})(\mathbb A_r+{\rm id})^{-1})^{-1},$$
and the Neumann series of the operator $({\rm id}+(p^{\alpha}-{\rm id})(\mathbb A_r+{\rm id})^{-1})^{-1}$ for $p$ sufficiently close to $1$.
\end{proof}

Now we can state two important asymptotic expansions on the Laplace transform $\p_\nu \widehat{U}(p,x_0)$ of the normal derivative $\partial_\nu U(t,x_0)$.
\begin{theorem}\label{theorem:flux:lowfreq:asympt:general}
Let Assumptions \ref{ass:coef}, \ref{ass:alpha:basic} and \ref{ass:g} hold.
Let $U$ be the solution of problem \eqref{eq:2D3D:ibvp}.
Fix any $x_0\in \partial\Omega$ and $N\in\mathbb{N}$.
Then there exist $C,p_0>0$ and $R(p)$ satisfying $|R(p)|\le C p^{N\underline{\alpha}}$ for all $p\in(0,p_0)$ and
\begin{equation}\label{eq:flux:lowfreq:asympt}
\p_\nu \widehat{U}(p,x_0)=p^{-M-1}\left(\sum_{k=3}^M \sum_{\ell=0}^{N-1} \p_\nu  u_{k,\ell}(p,x_0) + R(p)\right),\quad\forall p\in(0,p_0),
\end{equation}
where the functions $u_{k,0}(p,\cdot)\in W^{2,r}(\Om)\hookrightarrow C^{1,r_0}(\overline{\Om})$ for $k=3,\dots,M$ are the unique solutions of
\begin{equation}\label{eq:bvp:u0}
\left\{
\begin{aligned}
-\nabla\cdot(\sigma\nabla u_{k,0})+qu_{k,0}&=0,\quad\mbox{in }\Om,\\
u_{k,0}&=k!p^{M-k}\varphi_k,\quad\mbox{on }\partial\Omega,
\end{aligned}
\right.
\end{equation}
and $u_{k,\ell}(p,\cdot)\in W^{2,r}(\Om)\hookrightarrow C^{1,r_0}(\overline{\Om})$ for $k=3,\dots,M$ and $\ell\in\mathbb{N}$ are recursively defined to be the unique solutions of
\begin{equation}\label{eq:bvp:ui}
\left\{
\begin{aligned}
-\nabla\cdot(\sigma\nabla u_{k,\ell})+qu_{k,\ell}&=-p^{\alpha}\rho u_{k,\ell-1},\quad\mbox{in }\Om,\\
u_{k,\ell}&=0,\quad\mbox{on }\partial\Omega.
\end{aligned}
\right.
\end{equation}
\end{theorem}
\begin{proof}
Under Assumption \ref{ass:g}, we have
$$\widehat{g}(p,x)=\sum_{k=3}^M k! p^{-k-1} \varphi_k(x),\quad \forall (p,x)\in(0,\infty)\times\partial\Omega.$$
Hence, the identity \eqref{eq:normalder:Uhat} in Lemma \ref{lemma:representation:resolvents} gives
\begin{equation}\label{eq:lemma3.2:copy}
    \p_\nu \widehat{U}(p,x)=
\sum_{k=3}^M \p_\nu \left[\left({\rm id}-(\mathbb A_r+p^{\alpha})^{-1}p^{\alpha}\right)k! p^{-k-1}\mathcal{S} \varphi_k\right](x),\quad\forall (p,x)\in(0,\infty)\times\partial\Omega.
\end{equation}
Now fix any $N\in\mathbb{N}$ with $N\ge2$.
Upon letting
$$\mathbb{B}_N:=(\mathbb A_r + p^\alpha)^{-1}-\mathbb A_r^{-1}\textstyle\sum_{i=0}^{N-2}(-p^{\alpha}\mathbb A_r^{-1})^i,$$ we have the following identity
\begin{align}
{\rm id}-(\mathbb A_r + p^\alpha)^{-1}p^\alpha &= {\rm id}-\left(\mathbb A_r^{-1}\sum_{i=0}^{N-2}(-p^{\alpha}\mathbb A_r^{-1})^i +\mathbb B_N\right)p^\alpha \nonumber\\
&= \sum_{\ell=0}^{N-1}(-\mathbb A_r^{-1}p^{\alpha})^\ell-\mathbb B_N p^\alpha,\quad\forall p>0.\label{eq:lemma3.3:copy}
\end{align}
It follows from Lemma \ref{lemma:neumannseries} that there exist constants $C,p_0>0$ independent of $p$ satisfying
\begin{equation}\label{eqn:bound-Bn}
\|\mathbb B_N\|_{L^r(\Om)\to W^{2,r}(\Om)}\le C p^{(N-1)\underline{\alpha}},\quad \forall p\in(0,p_0).
\end{equation}
Meanwhile, from the definitions of $\mathcal{S}$ and $\mathbb A_r$, we deduce
\begin{align*}
u_{k,0}&=k! p^{M-k}\mathcal{S}\varphi_k,\\
u_{k,\ell} &= (-\mathbb A_r^{-1} p^{\alpha})u_{k,\ell-1},\quad \forall k\in\{3,\dots,M\},\, \forall\ell\in\mathbb{N}.
\end{align*}
Consequently, for all $k=3,\dots,M$ and $\ell\in\mathbb{N}\cup\{0\}$,
\begin{equation}\label{eq:intermsof:ukl}
\p_\nu (-\mathbb A_r^{-1}p^{\alpha})^\ell (k! p^{-k-1}\mathcal{S} \varphi_k)(x)=p^{-M-1}\p_\nu  u_{k,\ell}(p,x),\quad\forall (p,x)\in(0,\infty)\times\partial\Omega.
\end{equation}
By plugging \eqref{eq:lemma3.3:copy} and \eqref{eq:intermsof:ukl} into \eqref{eq:lemma3.2:copy}, we obtain the desired assertion \eqref{eq:flux:lowfreq:asympt} with the remainder
$$R(p)=-p^{M+1}\sum_{k=3}^M \p_\nu \left[\mathbb{B}_N p^\alpha (k!p^{-k-1}\mathcal{S}\varphi_k)\right](x_0),\quad\forall p\in(0,p_0).$$
From the estimate \eqref{eqn:bound-Bn}, we deduce
$$\|\mathbb B_N p^\alpha\|_{L^r(\Om)\to W^{2,r}(\Om)}\le C p^{N\underline{\alpha}},\quad \forall p\in(0,p_0).$$
This and the trace lemma imply
$$\|\p_\nu \left[\mathbb{B}_N p^\alpha (k!p^{-k-1}\mathcal{S}\varphi_k)\right]\|_{W^{1-1/r,r}(\p\Om)}\le C p^{N\underline{\alpha}-k-1},\quad\forall p\in(0,p_0),\,\forall k\in\{3,\dots,M\}.$$
In view of the Sobolev embedding $W^{1-1/r,r}(\p\Om)\hookrightarrow C^{0,r_0}(\p\Om)$, we deduce $$|R(p)|\le C p^{\min\{N\underline{\alpha}-k+M\,:\,3\le k\le M\}}= C p^{N\underline{\alpha}},\quad\forall p\in(0,p_0),$$
and hence the theorem for every $N\ge2$.
Also the assertion for $N=2$ directly implies that for $N=1$. This completes the proof of the theorem.
\end{proof}

\begin{remark}
\label{remark:Neq2:asympt}
Actually, we have $|\p_\nu  u_{k,\ell}(p,x_0)|=O(p^{M-k+\ell\underline{\alpha}})$ for all $p\in(0,p_0)$ from the expression \eqref{eq:intermsof:ukl} in the proof of Theorem \ref{theorem:flux:lowfreq:asympt:general}. Thus, for all $p\in(0,p_0)$, there holds
\begin{align*} \left|\sum_{k<M} \sum_{\ell\in\{0,1\}} \p_\nu u_{k,\ell}(p,x_0)\right|&=O\left(\sum_{k<M}\sum_{\ell\in\{0,1\}} p^{M-k+\ell\underline{\alpha}}\right)\\
&=O\left(p^{\min\{M-k+\ell\underline{\alpha}\,:\,k<M,\,\ell\in\{0,1\}\}}\right)=O(p).
\end{align*}
Hence the expansion \eqref{eq:flux:lowfreq:asympt} with $N=2$ implies
\begin{equation}\label{eq:main:asymptotics:for:uniqueness}
\p_\nu \widehat{U}(p,x_0)=p^{-M-1}\left( \sum_{\ell\in\{0,1\}} \p_\nu  u_{M,\ell}(p,x_0) + O(p+p^{2\underline{\alpha}})\right),\quad\forall p\in(0,p_0),
\end{equation}
which plays a crucial role in the proof of uniqueness of the inverse problem.
\end{remark}

The next result gives the asymptotic expansion of $\partial_\nu \widehat U(p,x_0)$ at $p=1$.

\begin{theorem}\label{theorem:flux:freqisone:asympt:general}
Let Assumptions \ref{ass:coef}, \ref{ass:alpha:basic} and \ref{ass:g} hold.
Let $U$ be the solution of problem \eqref{eq:2D3D:ibvp}.
Fix any $x_0\in \partial\Omega$ and $N\in\mathbb{N}$.
Then there exist $C,\varepsilon>0$ and $R(p)$ satisfying $|R(p)|\le C \|p^{\alpha}-{\rm id}\|_{L^\infty(\Om)}^{N}$ for all $p\in(1-\varepsilon,1+\varepsilon)$ and
\begin{equation}\label{eq:flux:lowfreq:asympt:pnear1}
\p_\nu \widehat{U}(p,x_0)=p^{-M-1}\left(\sum_{k=3}^M \sum_{\ell=0}^{N-1} \p_\nu  v_{k,\ell}(p,x_0) + R(p)\right),\quad\forall p\in(1-\varepsilon,1+\varepsilon),
\end{equation}
where $v_{k,0}(p,\cdot)\in W^{2,r}(\Om)\hookrightarrow C^{1,r_0}(\overline{\Om})$ for $k=3,\dots,M$ are the unique solutions of
\begin{equation}\label{eq:bvp:v0}
\left\{
\begin{aligned}
-\nabla\cdot(\sigma\nabla v_{k,0})+(q+\rho)v_{k,0}&=0,\quad\mbox{in }\Om,\\
v_{k,0}&=k!p^{M-k}\varphi_k,\quad\mbox{on }\partial\Omega,
\end{aligned}
\right.
\end{equation}
and $v_{k,\ell}(p,\cdot)\in W^{2,r}(\Om)\hookrightarrow C^{1,r_0}(\overline{\Om})$ for $k=3,\dots,M$ and $\ell\in\mathbb{N}$ are recursively defined to be the unique solutions of
\begin{equation}\label{eq:bvp:vi}
\left\{
\begin{aligned}
-\nabla\cdot(\sigma\nabla v_{k,\ell})+(q+\rho)v_{k,\ell}&=(1-p^{\alpha})\rho v_{k,\ell-1}(p,\cdot),\quad\mbox{in }\Om,\\
v_{k,\ell}&=0,\quad\mbox{on }\partial\Omega.
\end{aligned}
\right.
\end{equation}
\end{theorem}
\begin{proof}
By combining Lemmas \ref{lemma:representation:resolvents} and \ref{lemma:neumannseries:near:1} with an argument similar to that used for Theorem \ref{theorem:flux:lowfreq:asympt:general}, we obtain the desired assertion \eqref{eq:flux:lowfreq:asympt:pnear1}.
\end{proof}

\section{Proof of main results}\label{sec:inverse}
In this section, we present the detailed proofs of the main results stated in Section \ref{sec:main}. Throughout this section, $\mathbb A_r^i$, $i=1,2$, denotes the operator $\mathbb A_r$ associated with the coefficients $\sigma=\sigma^i$, $\rho=\rho^i$ and $q=q^i$, and $\mathcal{S}^i$ the corresponding solution operator, cf. \eqref{eq:S:bvp:wholedomain}.
The next result follows directly from the strong maximum principle and Hopf's lemma, and is used repeatedly below.
\begin{lemma}\label{lemma:maxprin:and:Hopf}
    Suppose that Assumption \ref{ass:coef} holds for $\sigma=\sigma^i$, $\rho=\rho^i$ and $q=q^i$ with $i=1,2$.
    If $f\in L^r(\Om)$ satisfies $f\ge0$ almost everywhere and $f$ does not uniformly vanish in the domain $\Omega$, then $\partial_\nu(\mathbb{A}_r^i)^{-1}f(x)<0$ for all $x\in\p\Om$ and $i=1,2$.
\end{lemma}
\begin{proof}
    Under Assumption \ref{ass:coef}, Lemma \ref{lemma:iso:Apluspalpha} with $p=0$ gives $(\mathbb{A}_r^i)^{-1}f\in W^{2,r}(\Om)\cap H_0^1(\Omega)$ for $i=1,2$.
    Thanks to the Sobolev embedding $W^{2,r}(\Om)\hookrightarrow C^{1,r_0}(\overline{\Om})$, with $0<r_0\leq 1-\frac{d}{r}$, the normal derivative $\p_\nu(\mathbb{A}_r^i)^{-1}f$ is well defined as a continuous function on $\p\Om$.
 Since $f\ge0$ everywhere and $f\not\equiv0$, by the strong maximum principle \cite[Theorems 3.5 and 9.6]{GT}, we have $(\mathbb{A}_r^i)^{-1}f(x)>0$ for all $x\in\Om$. Then the desired result follows from Hopf's lemma \cite[Lemma 3.4]{GT}.
\end{proof}

\subsection{Proof of Theorem \ref{theorem:flux:lowfreq:asympt:onept}}

In this part, we prove Theorem \ref{theorem:flux:lowfreq:asympt:onept}. The proof relies crucially on the asymptotic expansion in Theorem \ref{theorem:flux:lowfreq:asympt:general} with $N=2$ and the nonvanishing property of the normal derivative $\partial_\nu(\mathbb{A}_r^i)^{-1}f$ in Lemma \ref{lemma:maxprin:and:Hopf}.
\begin{proof}
Theorem \ref{theorem:flux:lowfreq:asympt:general} for $N=2$ and $U=U^i$ with $i=1,2$ gives, as stated in Remark \ref{remark:Neq2:asympt},
\begin{equation}\label{eq:pf:uniq:asympt-0}
\p_\nu \widehat{U}^i(p,x_0)=p^{-M^i-1}\left( \sum_{\ell=0}^{1} \p_\nu  u_{M^i,\ell}(p,x_0) + O(p+p^{2\underline{\alpha^i}})\right),\quad\mbox{as }p\to0^+.	
\end{equation}
For $i=1,2$, by the piecewise constancy of $\alpha_i$, there holds
$$\p_\nu  u_{M^i,1}^i=-\p_\nu  (\mathbb A_r^i)^{-1}[p^{\alpha^i} u_{M^i,0}^i]=-\sum_{j=0}^{n^i}p^{\alpha_j^i}\p_\nu  (\mathbb A_r^i)^{-1}\left[ u_{M^i,0}^i \chi_{\Om_j^i}\right].$$
By Theorem \ref{theorem:analyticity}, the map $t\mapsto \partial_\nu U^i(t,x_0)$ is analytic in $(0,\infty)$. Since the sequence $\{t_k\}_{k=1}^\infty$ consists of distinct numbers and accumulates to a positive number, the condition   $\sigma^1(x_0)\p_\nu  U^1(t_k,x_0)=\sigma^2(x_0)\p_\nu  U^2(t_k,x_0)$, $k\in\mathbb N$, from \eqref{t1a} and analytic continuation in time yield
$$\sigma^1(x_0)\p_\nu  U^1(t,x_0)=\sigma^2(x_0)\p_\nu  U^2(t,x_0),\quad \forall t>0.$$
Thus, by applying Laplace transform in time and using \eqref{eq:pf:uniq:asympt-0}, we obtain that as $p\to0^+$,
\begin{align}
&\hskip5mm \sigma^1(x_0){p^{-M^1-1}}\left(\p_\nu  u_{M^1,0}^1(x_0)-\sum_{j=0}^{n^1}p^{\alpha_j^1}\p_\nu  (\mathbb A_r^1)^{-1}\left[ u_{M^1,0}^1 \chi_{\Om_j^1}\right](x_0)+O(p+p^{2\alpha_{0}^1})\right)\nonumber \\
&=\sigma^2(x_0){p^{-M^2-1}}\left(\p_\nu  u_{M^2,0}^2(x_0)-\sum_{j=0}^{n^2}p^{\alpha_j^2}\p_\nu  (\mathbb A_r^2)^{-1}\left[ u_{M^2,0}^2 \chi_{\Om_j^2}\right](x_0)+O(p+p^{2\alpha_{0}^2})\right).\label{eq:thm32:mainid}
\end{align}
Next we prove the assertions in parts (a) and (b) separately.

\medskip
\noindent{\bf Part (a).}
The proof employs the relation \eqref{eq:thm32:mainid} essentially.
Without loss of generality, suppose that $M^1\ge M^2$ and $\varphi_{M^1}^1\ge0$ is not uniformly vanishing (recall that $\varphi_{M^1}^1$ is from the Dirichlet boundary data $g^1$).
Then we have $u_{M^1,0}^1(x)>0$ for all $x\in\Om$ by the strong maximum principle \cite[Theorem 3.5]{GT}, and thus, by Lemma \ref{lemma:maxprin:and:Hopf}, we conclude
\begin{equation}\label{eqn:sign-flux}
\p_\nu  (\mathbb A_r^1)^{-1}[ u_{M^1,0}^1 \chi_{\Om_j^1}](x_0)<0,\quad j=0,\ldots,n^1.
\end{equation}
First we prove $M^1=M^2$ by means of contradiction. Suppose to the contrary that $M^1> M^2$.
Then, by multiplying $p^{M^1+1}$ on both sides of the identity \eqref{eq:thm32:mainid} and then taking the limit $p\to 0^+$, we obtain $\p_\nu  u_{M^1,0}^1(p,x_0)=0$, since the function $u_{M^i,0}^i=M^i!\mathcal{S}^i\varphi^i_{M^i}$ is independent of $p$ for $i=1,2$.
Next, by multiplying $p^{M^1+1-\alpha_0^1}$ on both sides of \eqref{eq:thm32:mainid} and then taking the limit $p\to 0^+$, we obtain a contradiction: the left-hand side becomes $-\sigma^1(x_0)\p_\nu  (\mathbb{A}_r^1)^{-1}\left[ u_{M^1,0}^1 \chi_{\Om_0^1}\right](x_0)>0$, but the right-hand side becomes zero.
Therefore, we have $M^1=M^2$, which allows simplifying the identity \eqref{eq:thm32:mainid} into
$$\sum_{j=0}^{n^1}p^{\alpha_j^1}\p_\nu  (\mathbb A_r^1)^{-1}\left[ u_{M^1,0}^1 \chi_{\Om_j^1}\right](x_0)=\frac{\sigma^2(x_0)}{\sigma^1(x_0)}\sum_{j=0}^{n^2}p^{\alpha_j^2}\p_\nu  (\mathbb A_r^2)^{-1}\left[ u_{M^2,0}^2 \chi_{\Om_j^2}\right](x_0)+O(p^{\min\{1,2\alpha_{0}^1,2\alpha_{0}^2\}}).$$
Since $\{\alpha_j^i\}_{j=0}^{n^i}$ is strictly increasing for $i=1,2$ and by \eqref{eqn:sign-flux}, every coefficient of $p^{\alpha_j^i}$ on both sides is nonzero, one can easily show that $n^1=n^2$ and $\alpha_j^1=\alpha_j^2$ for all $j$, and hence part (a) of the theorem follows.

\medskip
\noindent{\bf Part (b).}
Let $\rho^1=\rho^2=:\rho$, $\sigma^1=\sigma^2=:\sigma$, $q^1=q^2=:q$ (so that $\mathbb A_r^1=\mathbb A_r^2=:\mathbb A_r$).
Also, let $M^1=M^2=:M$ and $\varphi_M^1=\varphi_M^2=:\varphi_M$.
From the proof of part (a), we deduce
\begin{equation*}
n^1=n^2=:n\quad \mbox{and}\quad \alpha_j^1=\alpha_j^2,\quad  j=0,1,\dots,n.
\end{equation*}
Consequently, $u_{M,0}^1=M!\mathcal{S}\varphi_{M}=u_{M,0}^2$.
Let $u_{M,0}^1=u_{M,0}^2=:u_{M,0}$.
By identifying the coefficients of the powers of $p$ on both sides of \eqref{eq:thm32:mainid}, we obtain
\begin{equation}\label{eq:coeff:identify}
\p_\nu  \mathbb A_r^{-1}\left[ u_{M,0} \chi_{\Om_j^1}\right](x_0)=\p_\nu  \mathbb A_r^{-1}\left[ u_{M,0} \chi_{\Om_j^2}\right](x_0), \quad j=0,1,\dots,n.
\end{equation}
Summing the identity \eqref{eq:coeff:identify} over $j\in\{0,1,\dots,k\}$ gives, for $k=0,1,\dots,n$,
\begin{equation}\label{eq:coeff:identify:cumulate}
\p_\nu  \mathbb A_r^{-1}\left[ u_{M,0} \chi_{\cup_{j=0}^k\Om_j^1}\right](x_0)=\p_\nu  \mathbb A_r^{-1}\left[ u_{M,0} \chi_{\cup_{j=0}^k\Om_j^2}\right](x_0).
\end{equation}
Meanwhile, using Lemma \ref{lemma:maxprin:and:Hopf} as in part (a), for all measurable subsets $A,B$ of $\Om$, the following strict monotonicity holds:
$$\mbox{if }A\subset B\mbox{ and }|B\backslash A|>0,\mbox{ then }0>\p_\nu  \mathbb A_r^{-1}\left[ u_{M,0} \chi_{A}\right](x_0)>\p_\nu  \mathbb A_r^{-1}\left[ u_{M,0} \chi_{B}\right](x_0).$$
The first assumption $\alpha^1\ge\alpha^2$ implies $\cup_{j=0}^k\Om_j^2\subset \cup_{j=0}^k\Om_j^1$ for all $k=0,1,\dots,n$, while the second assumption implies either $\Om_j^{1}\subset \Om_j^{2}$ or $\Om_j^{2}\subset \Om_j^{1}$ holds for all $j=0,1,\dots,n$.
By combining the strict monotonicity with the identities \eqref{eq:coeff:identify} and \eqref{eq:coeff:identify:cumulate}, we conclude $\Om_j^1=\Om_j^2$ for all $j=0,1,\dots,n$ if either the first or the second assumption holds true. This completes the proof of the theorem.
\end{proof}

\subsection{Proof of Theorem \ref{theorem:uniqueness:measfuns}}

The proof of the theorem proceeds by means of contradiction, and the main idea is to combine Hopf's lemma with the asymptotic expansion in Theorem \ref{theorem:flux:lowfreq:asympt:general}
and the continuity of Lebesgue measures.

\begin{proof}
Let $\mu$ be the Lebesgue measure in $\mathbb{R}^d$. From now on we assume that the condition \eqref{t1a} is fulfilled and we will prove that it implies the assertion \eqref{t1b} as well as the second claim of the theorem. Specifically, for  \eqref{t1b}, we prove that
$$\underline{\alpha^1}:=\operatorname{essinf}\{\alpha^1(x)\,:\,x\in\Om\}=\underline{\alpha^2}:=\operatorname{essinf}\{\alpha^2(x)\,:\,x\in\Om\}.$$
For the last statement of the theorem, under the additional assumptions, we shall prove that $\mu(\{x\in\Om\,:\,\alpha^1(x)\ne\alpha^2(x)\})=0$. The proof is divided into two steps, each dealing with one assertion.

\medskip
\noindent{\bf Step 1: proof of \eqref{t1b}.}
We prove $\underline{\alpha^1}=\underline{\alpha^2}$ by means of contradiction: Suppose to the contrary that $\underline{\alpha^1}\ne \underline{\alpha^2}$, say $\underline{\alpha^1}<\underline{\alpha^2}$.
Fix any $m\in(\underline{\alpha}^1,\min\{2\underline{\alpha^1},\underline{\alpha^2}\})$.
Then by the definitions of $\underline{\alpha^1}$ and $\underline{\alpha^2}$, we have $$\mu\left(\{x\in\Om\,:\, \alpha^1(x)<m\}\right)>0\quad\mbox{and}\quad\mu\left(\{x\in\Om\,:\, \alpha^2(x)<\underline{\alpha}^2\}\right)=0.$$
By the analytic continuation in time of the hypothesis as in the derivation of the identity  \eqref{eq:thm32:mainid}, there holds
$$\sigma^1(x_0)\p_\nu  \widehat{U}^1(p,x_0)=\sigma^2(x_0)\p_\nu  \widehat{U}^2(p,x_0)=:F(p),\quad \forall p>0.$$
Theorem \ref{theorem:flux:lowfreq:asympt:general} with $N=2$ leads to the following asymptotic relation (cf. Remark \ref{remark:Neq2:asympt}):
\begin{equation}\label{eq:F:twoasymptotics}
    F(p)= p^{-M^i-1}\sigma^i(x_0)\left(\p_\nu  u_{M^i,0}^i(x_0) - \p_\nu (\mathbb{A}_r^{i})^{-1}[p^{\alpha^i}u_{M^i,0}^i](x_0)+ O\big(p+p^{2\underline{\alpha^i}}\big)\right),
\end{equation}
as $p\to 0^+$ for $i=1,2$.
In particular, for $i=1,2$, we have
\begin{align}
\lim_{p\to 0^+} p^{M^i+1} |F(p)|&= \sigma^i(x_0)|\p_\nu  u_{M^i,0}^i(x_0)|<\infty,\label{eq:powerplus1}
\\
\lim_{p\to 0^+} p^{M^i} |F(p)|&= \sigma^i(x_0)\lim_{p\to 0^+} p^{-1}\left|\p_\nu  u_{M^i,0}^i(x_0) - \p_\nu (\mathbb{A}_r^{i})^{-1}[p^{\alpha^i}u_{M^i,0}^i](x_0)\right|=\infty.\label{eq:powerplus0}
\end{align}
The last relation of \eqref{eq:powerplus0} for the case  $\p_\nu  u_{M^i,0}^i(x_0)\ne0$ is immediate.
It suffices to show the last relation of \eqref{eq:powerplus0} for the case $\p_\nu  u_{M^i,0}^i(x_0)=0$. Indeed, since $\varphi_{M^i}^i$ is not sign-changing, $u_{M^i,0}^i$ is also not sign-changing. Then by Lemma \ref{lemma:maxprin:and:Hopf}, there holds for all $0<p<1$,
\begin{align*}
    |\p_\nu (\mathbb{A}_r^{i})^{-1}[p^{\alpha^i}u_{M^i,0}^i](x_0)|&\ge |\p_\nu (\mathbb{A}_r^{i})^{-1}[\chi_{\alpha^i<\frac{\underline{\alpha^i}+1}{2}}p^{\alpha^i}u_{M^i,0}^i](x_0)|\\
    &\ge p^{\frac{\underline{\alpha^i}+1}{2}}|\p_\nu (\mathbb{A}_r^{i})^{-1}[\chi_{\alpha^i<\frac{\underline{\alpha^i}+1}{2}}u_{M^i,0}^i](x_0)|,
\end{align*}
where the estimate $|\p_\nu (\mathbb{A}_r^{i})^{-1}[\chi_{\alpha^i<\frac{\underline{\alpha^i}+1}{2}}u_{M^i,0}^i](x_0)|>0$ holds due to the condition $\mu(\{x\in\Om\,:\,\alpha^i(x)<\frac{\underline{\alpha^i}+1}{2}\})>0$ and Lemma \ref{lemma:maxprin:and:Hopf} again. Thus the second assertion in \eqref{eq:powerplus0} holds also for the case $\p_\nu  u_{M^i,0}^i(x_0)=0$.
Since there can only be one integer $M^i$ satisfying both \eqref{eq:powerplus1} and \eqref{eq:powerplus0} simultaneously, we deduce $M^1=M^2=:M$.
Using this fact, the identity \eqref{eq:F:twoasymptotics} with $i=2$ and the choice $m<\underline{\alpha^2}$, we obtain
\begin{align}
&\lim_{p\to 0^+} p^{-m}\left(p^{M+1} F(p)-\sigma^2(x_0)\p_\nu  u_{M,0}^2(x_0)\right)\nonumber\\
= &\lim_{p\to0^+}\left(-\sigma^2(x_0)\p_\nu (\mathbb{A}_r^{2})^{-1}[p^{\alpha^2-m}u_{M^2,0}^2](x_0)+ O\big(p^{1-m}+p^{2\underline{\alpha^i}-m}\big)\right)=0.\label{eq:new:limit-m}
\end{align}
It follows from the relations $$\sigma^2(x_0)\p_\nu  u_{M^2,0}^2(x_0)=\lim_{p\to0^+}p^{M+1}F(p)=\sigma^1(x_0)\p_\nu  u_{M^1,0}^1(x_0),$$ and the identity \eqref{eq:F:twoasymptotics} with $i=1$ that
there holds, as $p\to 0^+$,
\begin{align*}
p^{M+1} F(p)-\sigma^2(x_0)\p_\nu  u_{M,0}^2(x_0)&=p^{M+1} F(p)-\sigma^1(x_0)\p_\nu  u_{M,0}^1(x_0)\\
&=-\sigma^1(x_0)\p_\nu (\mathbb{A}_r^{1})^{-1}[p^{\alpha^1}u_{M,0}^1](x_0)+ O\big(p+p^{2\underline{\alpha^1}}\big).
\end{align*}
By plugging this relation into the identity \eqref{eq:new:limit-m} and noting the choice $m<\min\{2\underline{\alpha^1},\underline{\alpha^2}\}\le\min\{2\underline{\alpha^1},1\}$, we arrive at
$$\lim_{p\to 0^+} p^{-m}\sigma^1(x_0)\p_\nu (\mathbb{A}_r^{1})^{-1}[p^{\alpha^1}u_{M,0}^1](x_0)=0.$$
This however contradicts the choice $m>\underline{\alpha^1}$ since Lemma \ref{lemma:maxprin:and:Hopf} and the condition $\mu(\{x\in\Om\,:\,\alpha^1(x)<m\})>0$ give the following relation: for all $0<p<1$,
$$p^{-m}|\p_\nu (\mathbb{A}_r^{1})^{-1}[p^{\alpha^1}u_{M,0}^1](x_0)|\ge|\p_\nu (\mathbb{A}_r^{1})^{-1}[\chi_{\alpha^1<m} u_{M,0}^1](x_0)|>0.$$
Thus, we conclude that $\underline{\alpha^1}=\underline{\alpha^2}$, and complete the proof of part (a).

\medskip
\noindent{\bf Step 2: Proof of part (b).}
We suppress the superscript $i$ as $\mathbb A_r^i=\mathbb A_r$ and $u_{M,0}^i=u_{M,0}$, $\underline{\alpha^i}=\underline{\alpha}$ etc., in view of the additional assumptions in the theorem.
We prove the assertion by contradiction.
Suppose to the contrary that $\mu(\{x\in\Om\,:\,\alpha^1(x)\ne\alpha^2(x)\})>0$ holds.
Then, by the continuity of measures, the ascending sequence of measurable sets $E_n:=\{x\in\Om\,:\,\alpha^1(x)-\alpha^2(x)>1/n\}$ satisfies
$$\lim_{n\to\infty}\mu(E_n)=\mu(\{x\in\Om\,:\,\alpha^1(x)-\alpha^2(x)\ne 0\})>0.$$
In other words, there exists some $N\in\mathbb{N}$ such that $\mu(E_n)\ge\mu(E_N)>0$ for all $n\ge N$.
Without loss of generality, we may assume $\varphi_M^1\le0$ on $\p\Om$.
Then we have $u_{M,0}<0$ in $\Om$ by the strong maximum principle.
Thus, from Lemma \ref{lemma:maxprin:and:Hopf}, we have
\begin{align*}
&\quad\p_\nu  \mathbb A_r^{-1}\left[(p^{\alpha^2}-p^{\alpha^1})u_{M,0}\right](x_0)\ge \p_\nu  \mathbb A_r^{-1}\left[\chi_{E_n}(p^{\alpha^2}-p^{\alpha^1})u_{M,0}\right](x_0)\\
&\ge (1-p^{1/n})\p_\nu  \mathbb A_r^{-1}\left[\chi_{E_n}p^{\alpha^2}u_{M,0}\right](x_0)\ge C'(1-p^{1/n})p^{\overline{\alpha^2}},\quad\forall n\ge N\mbox{ and }p\in(0,1),
\end{align*}
where $C' = \p_\nu  \mathbb A_r^{-1}\left[\chi_{E_N}u_{M,0}\right](x_0)>0 $.
Let $n_p=\lceil |\log p|\rceil$. Then, for all $n\ge n_p$, we have
\begin{equation*}
1-p^{1/n}\ge1-p^{1/|\log p|}=1-e^{-1},
\end{equation*}
which gives
$$\p_\nu  \mathbb A_r^{-1}\left[(p^{\alpha^2}-p^{\alpha^1})u_{M,0}\right](x_0)\ge C'(1-e^{-1})p^{\overline{\alpha^2}},\quad\forall p\in(0,1).$$
Meanwhile, the identity \eqref{eq:F:twoasymptotics} for $i=1$ and $i=2$ gives
$$|\p_\nu  \mathbb A_r^{-1}\left[(p^{\alpha^2}-p^{\alpha^1})u_{M,0}\right](x_0)|\le C (p+p^{2\underline{\alpha}}),\quad\forall p\in(0,1),$$
which contradicts the fact $\overline{\alpha^2}< \min\{1,2\underline{\alpha}\}$.
\end{proof}

\subsection{Proof of Theorem \ref{theorem:stability}}\label{ssec:thm:stability}
The proofs of Theorems \ref{theorem:stability} and \ref{theorem:stability:fullbdydata} rely essentially on an integral representation of the flux data in the Laplace domain given in Lemma \ref{lem:int-repres} below, and the $W^{2,r}(\Omega)$ regularity and the strong maximum principle for problem \eqref{eq:bvp:vtild:total:for:stabliliy} below.

First we give two preliminary results. The first result gives an important integral representation of the derivative of $\partial_\nu \widehat {U}(p,x_0)$ at $p=1$.
\begin{lemma}\label{lem:int-repres}
Let $U$ be the solution of problem \eqref{eq:2D3D:ibvp} under
Assumptions \ref{ass:coef}, \ref{ass:alpha:basic} and \ref{ass:g}. Then the following identity holds
\begin{equation}\label{eq:limit:diffp1}
    \lim_{p\to1}\frac{\p_\nu \widehat{U}(p,x_0)-\p_\nu \widehat{U}(1,x_0)}{p-1}=\int_{0}^\infty \p_\nu  U(t,x_0)(te^{-t})\d t.
\end{equation}
\end{lemma}
\begin{proof}
For all $p\in(0,2)\backslash\{1\}$ and $t>0$, Taylor's formula gives
$$|e^{-pt}-e^{-t}+(p-1)te^{-t}|\le \tfrac{1}{2}|p-1|^2 t^2 e^{-ct},$$
where the constant $c$ lies between $p$ and $1$.
Thus, by the definition of the Laplace transform, we have
\begin{align}
    &\hskip5mm \left|\p_\nu \widehat{U}(p,x_0)-\p_\nu \widehat{U}(1,x_0) + (p-1)\int_{0}^\infty \p_\nu  U(t,x_0){(te^{-t})}\d t\right|\nonumber\\
&\le \frac{|p-1|^2}{2}\int_0^\infty|\p_\nu  U(t,x_0)|t^2e^{-ct}\d t,\quad\forall p\in(0,2).\label{eq:CSineq:stability}
\end{align}
Dividing both sides of the identity \eqref{eq:CSineq:stability} by $p-1$ and then taking the limit $p\to1$ give the desired assertion.
\end{proof}

\begin{lemma}\label{lemma:linf:numb}
    Let $\alpha$ satisfy Assumption \ref{ass:alpha:basic}. Then there holds $$\|p^{\alpha}-1-(p-1)\alpha \|_{L^\infty(\Om)}\le 
    |p-1|^2,\quad\forall p\in(0.5,1.5).$$
\end{lemma}
\begin{proof}
For each $\beta\in(0,1)$ and $x\in(-1,1)$, Taylor expansion gives
$$(1+x)^{\beta}=1+\beta x+ \beta(\beta-1)(1+x^*)^{\beta-2} x^2,$$
for some $x^*$ between $0$ and $x$. Thus, for all $\beta\in(0,1)$ and $x\in(-0.5,0.5)$, there holds
\begin{equation*}
        |(1+x)^{\beta}-1-\beta x|\le 
        (0.5)^2(0.5)^{-2}x^2 = x^2.
    \end{equation*}
Using this inequality with $x=p-1$ and then taking supremum over $\beta\in(0,1)$ gives the desired inequality.
\end{proof}

Now we can state the proof of Theorem \ref{theorem:stability}.
\begin{proof}
Theorem \ref{theorem:flux:freqisone:asympt:general} with $N=2$ and $U=U^i$ with $i=1,2$ gives
\begin{equation}\label{eq:pf:uniq:asympt}
\p_\nu \widehat{U}^i(p,x_0)=p^{-M-1}\left( \sum_{k=3}^M\sum_{\ell=0}^{1} \p_\nu  v_{k,\ell}^i(p,x_0) + O(\|1-p^{\alpha}\|_{L^\infty(\Om)}^2)\right),\quad\mbox{as }p\to1,	
\end{equation}
with $v_{k,0}^i$ and $v_{k,1}^i$ satisfying \eqref{eq:bvp:v0} and \eqref{eq:bvp:vi}, respectively. Note that $v_{k,0}^i$ is independent of $\alpha^i$ and we write $v_{k,0}$ directly below by dropping the superscript $i$.
By Lemma \ref{lemma:linf:numb}, we may replace $O(\|1-p^{\alpha}\|_{L^\infty(\Om)}^2)$ by $O(|p-1|^2)$. Consequently, we obtain, as $p\to1$,
\begin{align}
&\hskip5mm \p_\nu \widehat{U}^1(p,x_0)-\p_\nu \widehat{U}^2(p,x_0) = p^{-M-1}\left(\sum_{k=3}^M\left(\p_\nu  v_{k,1}^1(p,x_0)-\p_\nu  v_{k,1}^2(p,x_0) \right)+ O(|p-1|^2)\right).\label{eq:thm32:mainid:stable}
\end{align}
For $k=3,\dots,M$ and $p>0$, the function $w_k:=v_{k,1}^1(p,\cdot)-v_{k,1}^2(p,\cdot)$ satisfies
\begin{equation}\nonumber
\left\{
\begin{aligned}
-\nabla\cdot(\sigma\nabla w_k)+(q+\rho)w_k&=-(p^{\alpha^1}-p^{\alpha^2})\rho v_{k,0}(p,\cdot),\quad\mbox{in }\Om,\\
w_k&=0,\quad\mbox{on }\partial\Omega.
\end{aligned}
\right.
\end{equation}
From Lemma \ref{lemma:iso:Apluspalpha} with $\kappa=r$ and $p=1$, we have
$$\|w_k-(1-p)\widetilde{v}_{k,1}\|_{W^{2,r}(\Om)}\le C\|(p^{\alpha^1}-p^{\alpha^2})\rho v_{k,0}(p,\cdot) + (1-p)(\alpha^1-\alpha^2)\rho v_{k,0}(1,\cdot)\|_{L^r(\Om)},\quad\forall p>0,$$
where the constant $C>0$ is independent of $p$, and
the function $\widetilde{v}_{k,1}$ satisfies
\begin{equation}\label{eq:bvp:vi:tilde}
\left\{
\begin{aligned}
-\nabla\cdot(\sigma\nabla \widetilde{v}_{k,1})+(q+\rho)\widetilde{v}_{k,1}&=(\alpha^1-\alpha^2)\rho v_{k,0}(1,\cdot),\quad\mbox{in }\Om,\\
\widetilde{v}_{k,1}&=0,\quad\mbox{on }\partial\Omega.
\end{aligned}
\right.
\end{equation}
Note that $v_{k,0}(p,x)-v_{k,0}(1,x)= (p^{M-k}-1)v_{k,0}(1,x)$ for all $p>0$, $x\in\Om$ and $k=3,\dots,M$. From this and Lemma \ref{lemma:linf:numb}, we obtain $$\|w_k-(1-p)\widetilde{v}_{k,1}\|_{W^{2,r}(\Om)}=O(|p-1|^2).$$
By the Sobolev embedding $W^{2,r}(\Om)\hookrightarrow C^{1,r_0}(\overline{\Om})$ with $0<r_0\le 1-\frac{d}{r}$,
we deduce
\begin{equation*}
\p_\nu  v_{k,1}^1(p,x_0)-\p_\nu  v_{k,1}^2(p,x_0)= (1-p)\p_\nu \widetilde{v}_{k,1}(x_0) + O(|p-1|^2),\quad\mbox{as }p\to1.
\end{equation*}
Moreover, we obtain $\sum_{k=3}^M\p_\nu \widetilde{v}_{k,1}(x_0)=\p_\nu  \widetilde{v}(x_0)$ for the solution $\widetilde{v}$ to
\begin{equation}\label{eq:bvp:v:total:for:stabliliy}
\left\{
\begin{aligned}
-\nabla\cdot(\sigma\nabla \widetilde{v})+(q+\rho)\widetilde{v}&=(\alpha^1-\alpha^2)\rho \sum_{k=3}^M v_{k,0}(1,\cdot),\quad\mbox{in }\Om,\\
\widetilde{v}&=0,\quad\mbox{on }\partial\Omega.
\end{aligned}
\right.
\end{equation}
By dividing both sides of \eqref{eq:thm32:mainid:stable} by $p-1$, taking the limit $p\to1$ and then comparing the identity with the identity \eqref{eq:limit:diffp1}, we arrive at
\begin{align} -\p_\nu  \widetilde{v}(x_0)&=\lim_{p\to1}\frac{\p_\nu \widehat{U}^1(p,x_0)-\p_\nu \widehat{U}^2(p,x_0)}{p-1}\nonumber\\
&=\int_{0}^\infty (\p_\nu  U^1(t,x_0)-\p_\nu  U^2(t,x_0))(te^{-t})\d t. \label{tata}
\end{align}
Fix $\{\alpha^{2,\ell}\}_{\ell\in\mathbb N}$ to be an $L^r(\Om)$-convergent sequence of  functions satisfying Assumption \ref{ass:alpha:basic} with $\alpha^{2,\ell}\geq \alpha^1$, $\ell\in\mathbb N$. For all $\ell\in\mathbb N$, let $U^{2,\ell}$
be the solution of problem \eqref{eq:2D3D:ibvp} with $\alpha=\alpha^{2,\ell}$. For each $\ell\in\mathbb N$, by choosing $\alpha^2=\alpha^{2,\ell}$, we  obtain from \eqref{tata} the following identity
\begin{equation} \label{eq:totakelimit}
-\p_\nu  \widetilde{v}_\ell(x_0)=\int_{0}^\infty (\p_\nu  U^1(t,x_0)-\p_\nu  U^{2,\ell}(t,x_0))(te^{-t})\d t,\quad \ell\in\mathbb N,
\end{equation}
where $\widetilde{v}_\ell$ solves problem \eqref{eq:bvp:v:total:for:stabliliy} with $\alpha^2=\alpha^{2,\ell}$.
To prove the desired assertion, we proceed by means of contradiction.  Suppose to the contrary that the right-hand side of \eqref{eq:totakelimit} converges to zero as $\ell$ tends to infinity.
We prove that $\|\alpha^{2,\ell}-\alpha^1\|_{L^r(\Om)}\to0$ as $\ell\to+\infty$ by means of contradiction. Suppose that $\|\alpha^{2,\ell}-\alpha^3\|_{L^r(\Om)}\to0$ as $\ell\to+\infty$  for some $\alpha^3\neq\alpha^1$ satisfying $\alpha^3\geq\alpha^1$. Let $\widetilde{v}_*$ be the solution of problem \eqref{eq:bvp:v:total:for:stabliliy} with $\alpha^2=\alpha^3$.
Then, by recalling the condition $r>d$, the continuous Sobolev embedding $W^{1-1/r,r}(\partial\Omega)\hookrightarrow C(\partial\Omega)$ and trace lemma, we have
\begin{equation}\label{eq:tracelemma:vtild}
    |\p_\nu  \widetilde{v}_\ell(x_0)-\p_\nu  \widetilde{v}_*(x_0)|\le C\|\p_\nu  \widetilde{v}_\ell-\p_\nu  \widetilde{v}_*\|_{W^{1-1/r,r}(\partial\Omega)}
    \leq C \|\widetilde{v}_\ell-\widetilde{v}_*\|_{W^{2,r}(\Om)},\quad \ell\in\mathbb N.
\end{equation}
Also, by applying Lemma \ref{lemma:iso:Apluspalpha} with $p=1$ to problem \eqref{eq:bvp:v:total:for:stabliliy} with $\alpha^1-\alpha^2$ replaced by $\alpha^3-\alpha^{2,\ell}$, we have
\begin{equation*}
    \|\widetilde{v}_\ell-\widetilde{v}_*\|_{W^{2,r}(\Om)} \leq C\|(\alpha^3-\alpha^{2,\ell})\rho v\|_{L^r(\Om)},\quad \ell\in\mathbb N.
\end{equation*}
Consequently, we deduce
\begin{equation*}
|\p_\nu  \widetilde{v}_\ell(x_0)-\p_\nu  \widetilde{v}_*(x_0)|\le  C\|(\alpha^3-\alpha^{2,\ell})\rho v\|_{L^r(\Om)},\quad \ell\in\mathbb N,
\end{equation*}
which tends to zero as $\ell\to\infty$. Meanwhile, taking the limit as $\ell\to\infty$ on both sides of the identity \eqref{eq:totakelimit} gives $\lim_{\ell\to\infty} \p_\nu  \widetilde{v}_\ell(x_0)=0$. Combining these estimates with the triangle inequality gives
\begin{equation*}
|\p_\nu  \widetilde{v}_*(x_0)|\leq \lim_{\ell\to \infty}\big(|\p_\nu  \widetilde{v}_\ell(x_0)-\p_\nu  \widetilde{v}_*(x_0)| + |\p_\nu  \widetilde{v}_\ell(x_0)|\big) =0.
\end{equation*}
Next we derive the contradiction.
Note that the function ${v}:=\sum_{k=3}^M v_{k,0}(1,\cdot)$ satisfies
\begin{equation}\label{eq:bvp:vtild:total:for:stabliliy}
\left\{
\begin{aligned}
-\nabla\cdot(\sigma\nabla {v})+(q+\rho){v}&=0,\quad\mbox{in }\Om,\\
{v}&= \sum_{k=3}^M k!\varphi_k,\quad\mbox{on }\partial\Omega.
\end{aligned}
\right.
\end{equation}
It follows from the hypothesis of the theorem that the boundary condition $\textstyle\sum_{k=3}^M k!\varphi_k$ is neither everywhere-vanishing nor sign-changing.
Then, by the strong maximum principle for problem \eqref{eq:bvp:vtild:total:for:stabliliy}, we have $v(x)>0$ for all $x\in\Om$.
Finally, recalling that $\alpha^3\geq\alpha^1$ and $\alpha^3\neq\alpha^1$, from Hopf's lemma for problem \eqref{eq:bvp:v:total:for:stabliliy} with $\alpha^2=\alpha^3$, we get $\p_\nu  \widetilde{v}_*(x_0)\ne0$, which contradicts the previous statement and completes the proof of the theorem.
\end{proof}

\subsection{Proof of Theorem \ref{theorem:stability:fullbdydata}}
We use the notation of Section \ref{ssec:thm:stability}. Let $\tilde{v}$ be the solution of \eqref{eq:bvp:v:total:for:stabliliy}. One can derive the following identity in the same way as in the proof of Theorem \ref{theorem:stability}:
\begin{equation} \label{eq:totakelimit:fullmeas}
-\p_\nu  \widetilde{v}(x)=\int_{0}^\infty (\p_\nu  U^1(t,x)-\p_\nu  U^{2}(t,x))(te^{-t})\,\d t,\quad x\in\partial\Omega,
\end{equation}
with $\widetilde v$ being the solution to \eqref{eq:bvp:v:total:for:stabliliy}.
Define the auxiliary function $w\in H^2(\Om)$ to be the  solution to
\begin{equation}\label{eq:bvp:w:aux}
\left\{
\begin{aligned}
-\nabla\cdot(\sigma\nabla w)+(q+\rho)w&=0,\quad\mbox{in }\Om,\\
w&=1,\quad\mbox{on }\partial\Omega.
\end{aligned}
\right.
\end{equation}
By the strong maximum principle, we have $w\geq 1$ in $\Omega$.
Then, by multiplying $w$ on both sides of equation \eqref{eq:bvp:v:total:for:stabliliy}, integrating over the domain $\Omega$ and integrating by parts, we obtain the identity
$$\int_\Om \sigma \nabla\widetilde{v}\cdot\nabla w \,{\rm d}x + \int_\Om (q+\rho)\widetilde{v}w \,{\rm d}x- \int_{\p\Om}\sigma (\p_\nu \widetilde{v})w \,{\rm d}s= \int_\Om (\alpha^1-\alpha^2)\rho v w \,{\rm d}x,$$
where $v$ is the solution to problem \eqref{eq:bvp:vtild:total:for:stabliliy}.
Using the identity \eqref{eq:bvp:w:aux}, we arrive at
$$\int_\Om (\alpha^1-\alpha^2)\rho v w\,{\rm d} x = - \int_{\p\Om}\sigma \p_\nu \widetilde{v}\, {\rm d }s.$$
By the strong maximum principle for problems \eqref{eq:bvp:vtild:total:for:stabliliy} and \eqref{eq:bvp:w:aux}, $v$ is not sign-changing and there hold $|v|\ge v_0>0$ and $w\ge1$ in $\Om$.  Using these two facts, the condition $\operatorname{essinf}\{\rho(x)\,:\,x\in\Omega\}>0$ from Assumption \ref{ass:coef} and the monotonicity condition $\alpha^1\ge\alpha^2$ gives
$$\|\alpha^1-\alpha^2\|_{L^1(\Om)}\leq C\|\p_\nu \widetilde{v}\|_{L^1(\p\Omega)},$$
with $C>0$ depending on $\sigma$, $\rho$, $\Omega$, $q$ and  $g$.
This completes the proof of the theorem.
\qed

\section{Conclusion}\label{sec:concl}
In this work we have established several new uniqueness, (semi-)continuity and stability results for recovering a spatially variable order function in the standard variable-order subdiffusion model from the flux observations on the boundary, either at one single point or over the whole boundary, for either piecewise constant order or general variable order. These uniqueness results substantially extend existing results in both one-dimensional case \cite{HJK:2025:ISDVO} and multi-dimensional case \cite{IkehataKian:2023,Kian:2018:TFD}, and to the best of our knowledge, the continuity and stability results are new in the context. Theoretically one important issue is to relax the monotonicity assumption on the variable order, which has played a vital role in the analysis. Numerically, one outstanding issue is to develop efficient algorithms for recovering the variable order from the boundary data and to establish the convergence of the algorithms, which so far have not been attempted.

\bibliographystyle{abbrv}
\bibliography{reference2}
\end{document}